\documentclass[12pt]{amsart}
\usepackage[T1]{fontenc}
\usepackage{amsmath, amssymb, latexsym, amsthm, mathrsfs, color, mathtools, tikz,tikz-cd, graphicx, caption, stackengine,url,accents, scrextend, bm, adjustbox}
\usepackage[symbol]{footmisc}
\usepackage[top=2.5cm, bottom=2.5cm, left=2.5cm, right=2.5cm]{geometry}
\usepackage{csquotes}

\stackMath
\newcommand\tsup[2][2]{%
 \def\useanchorwidth{T}%
  \ifnum#1>1%
    \stackon[-1pt]{\tsup[\numexpr#1-1\relax]{#2}}{\hspace{1pt}\scriptstyle\sim}%
  \else%
    \stackon[.5pt]{#2}{\hspace{1pt}\scriptstyle\sim}%
  \fi%
}

\newcommand{\nc}{\newcommand}
\nc{\Opn}{\mathrm{O}}

\DeclareMathOperator{\Int}{Int}
\DeclareMathOperator{\Intc}{Int_c}

\nc{\cO}{\mathcal{O}}

\newcommand{\sgg}{{\mathsf{S}_1(\Ga,\Ga)}}

\newcommand{\sgo}{{\mathsf{S}_1(\Ga,\Op)}}

\newcommand{\sgw}{{\mathsf{S}_1(\Ga,\Om)}}
\newcommand{\sww}{\mathsf{S}_1(\Omega,\Omega)}
\newcommand{\sfinww}{\sfin(\Omega,\Omega)}

\nc{\ram}{\mathsf{Ramsey}}
\nc{\soo}[1]{\mathsf{S}_1(\Op,\Op)}
\nc{\swl}[1]{\mathsf{S}_1(\Om,\Lambda)}
\nc{\swg}[1]{\mathsf{S}_1(\Om,\Ga)}
\nc{\goo}[1]{\gone(\Op,\Op)}
\nc{\gwo}{\gone(\Om,\cO)}
\nc{\gwl}[1]{\gone(\Om(#1),\Lambda(#1))}
\nc{\soox}[1]{\mathsf{S}_1(\cO(#1),\cO(#1))}
\nc{\swlx}[1]{\mathsf{S}_1(\Om(#1),\Lambda(#1))}
\nc{\goox}[1]{\gone(\cO(#1),\cO(#1))}
\nc{\gwox}[1]{\gone(\Om(#1),\cO(#1))}
\nc{\gwlx}[1]{\gone(\Om(#1),\Lambda(#1))}
\nc{\mich}[1]{{(#1)_\mathrm{M}}}
\nc{\michk}[2]{{(#1)^{#2}_\mathrm{M}}}
\nc{\PNM}{{\PN_\mathrm{M}}}
\nc{\sfinwl}{\mathsf{S}_{\mathrm{fin}}(\Om,\Lambda)}
\nc{\gfinoo}{\gfin(\cO,\cO)}
\nc{\gfinwo}{\gfin(\Om,\cO)}
\nc{\gfinwl}{\gfin(\Om,\Lambda)}
\nc{\sfinoox}[1]{\mathsf{S}_{\mathrm{fin}}(\cO(#1),\cO(#1))}
\nc{\sfinwlx}[1]{\mathsf{S}_{\mathrm{fin}}(\Om,\Lambda)}
\nc{\sfinwwx}[1]{\mathsf{S}_{\mathrm{fin}}(\Om(#1),\Om(#1))}
\nc{\gfinoox}[1]{\gfin(\cO(#1),\cO(#1))}
\nc{\gfinwox}[1]{\gfin(\Om(#1),\cO(#1))}
\nc{\gfinwlx}[1]{\gfin(\Om(#1),\Lambda(#1))}

\DeclareMathOperator{\fin}{Fin}

\nc{\mc}{\mathcal}
\nc{\thusfar}{\my{--- Edited thus far ---}}
\nc{\Omctbl}{\Omega_\mathrm{ctbl}}
\nc{\sctblww}{{\mathsf{S}_1(\Omctbl,\Om)}}
\nc{\lei}{\le^\oo}
\nc{\sqsubs}{\sqsubseteq^*}
\nc{\card}[1]{\left|#1\right|}
\nc{\medcard}[1]{\biggl|\,#1\,\biggr|}
\nc{\smallcard}[1]{|\,#1\,|}
\nc{\bds}{bidirectional $\roth$-scale}
\nc{\bfP}{\mathbf{P}}
\nc{\bfQ}{\mathbf{Q}}
\nc{\bbT}{\mathbb{T}}
\nc{\bbZ}{\mathbb{Z}}
\nc{\bbN}{\mathbb{N}}
\nc{\bbC}{\mathbb{C}}
\nc{\beq}{\begin{equation}}
\nc{\eeq}{\end{equation}}
\nc{\beqs}{\begin{equation*}}
\nc{\eeqs}{\end{equation*}}
\nc{\mbq}{\mb{?}}
\nc{\mb}[1]{{\mbox{\textbf{#1}}}}
\nc{\nop}{$\times$}
\nc{\fbn}{\!\!\fbox{\!\nop\!}\!\!}
\nc{\yup}{\checkmark}
\nc{\forces}{\Vdash}
\nc{\name}[1]{\dot{#1}}
\nc{\tf}{\my{FINISHED THUS FAR}}
\nc{\FU}{Fr\'echet--Urysohn}
\nc{\gs}{$\gamma$~space}
\nc{\Gab}{\Gamma_{\mathrm{B}}}
\nc{\Omb}{\Omega_{\mathrm{B}}}
\nc{\Ga}{\Gamma}
\nc{\Om}{\Omega}
\nc{\smallbinom}[2]{\begin{psmallmatrix} #1\\ #2 \end{psmallmatrix}}
\nc{\bgamma}{\smallbinom{\Om}{\Ga}}
\newcommand{\two}{\{0,1\}}
\nc{\productive}[2]{(#1,\allowbreak #2)^\x}
\nc{\prdct}[1]{#1^\x}
\nc{\Sel}{\mathsf{S}}
\nc{\sset}[2]{\{\,#1 : #2\,\}}
\nc{\smb}[1]{{\!\!\mb{#1}\!\!}}
\nc{\medset}[2]{{\biggl\{\,#1 : #2\,\biggr\}}}
\nc{\smallmedset}[2]{{\bigl\{\,#1 : #2\,\bigr\}}}
\nc{\set}[2]{{\left\{\,#1 : #2\,\right\}}}
\nc{\seq}[2]{{\la\, #1 : #2\,\ra}}
\nc{\eseq}[1]{#1_0, \allowbreak #1_1, \allowbreak\dotsc} 
\nc{\CH}{\mathsf{CH}}
\nc{\NCF}{\mathsf{NCF}}
\nc{\ZFC}{\mathsf{ZFC}}

\nc{\eseqint}[3]{#1_{#2}, \allowbreak\dotsc,\allowbreak #1_{#3}} 
\nc{\eseqstart}[2]{#1_{#2},\allowbreak #1_{#2+1},\dotsc } 
\nc{\eprod}[1]{#1_{1}\times \allowbreak#1_{2}\times\dotsb}

\nc{\shortprod}[1]{\prod_{n=1}^\infty{#1}_n}

\nc{\eprodint}[3]{#1_{#2}\times \allowbreak\dotsb\times\allowbreak #1_{#3}}

\nc{\seleseq}[1]{#1_1\in \mathcal{#1}_1, \allowbreak #1_2\in \mathcal{#1}_2, \allowbreak\dotsc}
\nc{\cube}{(\Cantor)^\bbN}
\nc{\Match}{\op{Match}}
\nc{\concat}[1]{\hat{\phantom{a}}\langle #1\rangle}
\nc{\poset}{\mathbb{P}}
\nc{\fn}[1]{{\op{Fn}(#1\times\w,2)}}
\nc{\linadd}{\op{linadd}}
\nc{\nonprod}{\non^\x}
\nc{\alephes}{{\aleph_0}}
\nc{\my}[1]{{\color{red}{#1}}}
\nc{\later}[1]{{\color{green} #1}}
\nc{\BTs}[1]{{\color{green} #1 (BT)}}
\nc{\Cp}{\op{C}_\mathrm{p}}
\nc{\Bp}{\op{B}_p}
\nc{\Pa}[8]{\bibitem{#1} {#2}, \emph{#3}, {#4} \textbf{#5} ({#6}), {#7}--{#8}.}
\nc{\tPa}[5]{\bibitem{#1} {#2}, \emph{#3}, {#4}, to appear.}
\nc{\sPa}[4]{\bibitem{#1} {#2}, \emph{#3}, {#4}, submitted.}
\nc{\Bc}[9]{\bibitem{#1} {#2}, \emph{#3}, in: \textbf{#4} (#5), #6 #7, #8--#9.}
\nc{\fD}{\mathfrak{D}}
\nc{\fX}{\mathfrak{X}}
\nc{\Onbd}{\Op_{\mathrm{nbd}}} 
\nc{\Omnb}{\Om_{\mathrm{nbd}}} 
\nc{\od}{\mathfrak{od}}
\nc{\Setting}[7]{\xymatrix@R=4pt@C=7pt{#1\ar@{-}[r]&#2\ar@{-}[r]&#3\\&#4\ar@{-}[u]\\
#5\ar@{-}[uu]\ar@{-}[r] & #6\ar@{-}[u]\ar@{-}[r] & #7\ar@{-}[uu]}}
\nc{\mx}[1]{\begin{matrix}#1\end{matrix}}
\nc{\plim}{p\txt{-}\lim}
\nc{\Bgp}{{\Z^\bbN}}
\nc{\Cgp}{{{\Z_2}^\bbN}}
\nc{\Cite}[1]{\textbf{[#1]}}
\nc{\Next}[1]{{#1^+}}
\nc{\cFin}{\mathrm{cF}}
\nc{\scsp}{\text{-scale space}}
\nc{\cfn}{\text{cofinal}\ }
\nc{\Con}{\text{Concentrated}}
\nc{\Lind}{\text{Lindel\"of}\,}
\nc{\con}{\text{-Concentrated}}
\nc{\lind}{\text{-Lindel\"of}\,}
\nc{\ctbl}{\text{countably }\allowbreak}
\nc{\Hur}{\text{Hurewicz}}
\nc{\intvl}[2]{{[#1(#2),\allowbreak #1(#2\!+\!1))}}
\nc{\Bdd}{\mathbf{B}}
\nc{\Dfin}{\mathfrak{D}_\mathrm{fin}}
\nc{\grbl}{{\mbox{\textit{\tiny gp}}}}
\nc{\bbP}{\mathbb{P}}
\nc{\BOfat}{\B_{\Om_{\mathrm{fat}}}}
\nc{\Bgood}{\B_{\mathrm{good}}}
\nc{\compactN}{\cl{\mathbb{N}}}
\nc{\blocks}[2]{\op{cl}_{#2}(#1)}
\nc{\blocksplus}[2]{\op{cl}^+_{#2}(#1)}
\nc{\arx}[1]{\texttt{http://arxiv.org/math/#1}}
\nc{\bq}{\begin{quote}}
\nc{\eq}{\end{quote}}
\nc{\cl}[1]{\overline{#1}}
\nc{\Cl}[2]{\mathrm{cl}_{#1}(#2)}
\nc{\MA}{Martin's Axiom}
\nc{\Bfat}{\B_\mathrm{fat}}
\nc{\inv}{^{-1}}
\nc{\Cantor}{{\two^\bbN}}
\nc{\bP}{\mathbf{P}}
\nc{\bof}{\op{\fb}}
\nc{\dof}{\op{\fd}}
\nc{\bofF}{\bof(\cF)}
\nc{\sr}[3]{\underset{\mbox{#3}}{\mbox{#1}}}
\nc{\gp}{\binom{\Om}{\Ga}}
\nc{\gpsmall}{\mbox{$\gp$}}
\nc{\gig}{\gimel}
\nc{\gns}{\sone(\Om,\gig)}
\nc{\nsr}[2]{#1}
\nc{\Srg}{{\mathbb{S}}}
\nc{\Srgs}{{\mathbb{S}^*}}
\nc{\NN}{{\w^{\w}}}
\nc{\ZN}{{\Z^{\bbN}}}
\nc{\NNup}{{\bbN^{\uparrow\bbN}}}
\nc{\NNupb}{{b^{\uparrow\bbN}}}
\nc{\Pof}{\op{P}}
\nc{\PN}{{\Pof(\w)}}
\nc{\rothx}[1]{{[#1]^{\mbox{\tiny $\infty$}}}}
\nc{\tx}{{\tilde{x}}}
\nc{\roth}{{[\w]^{\w}}}
\nc{\roths}{{[b]^{\mbox{\tiny $\infty$}}}} 
\nc{\Fin}{\mathrm{Fin}}
\nc{\ici}{[\bbN]^{ \infty, \infty}}
\nc{\Inc}{{\compactN^{\uparrow\bbN}}}
\nc{\powInc}[1]{{\big(\Inc\big)^{#1}}}
\nc{\powFin}[1]{{\big(\Fin\big)^{#1}}}
\nc{\powPN}[1]{{\big(\PN\big)^{#1}}}
\nc{\NcompactN}{{\compactN^\bbN}}
\nc{\Uarrow}{\smash{\big\uparrow}}
\nc{\LE}{\preccurlyeq}
\nc{\GE}{\succcurlyeq}
\nc{\op}{\operatorname}
\nc{\im}{\op{Im}}
\nc{\Span}{\op{span}}
\nc{\maxfin}{\op{maxfin}}
\nc{\ran}{\op{range}}
\nc{\iso}{\cong}
\nc{\Madd}{{\M}^\star}
\nc{\cI}{\mathcal{I}}
\nc{\cJ}{\mathcal{J}}
\nc{\scrA}{\mathscr{A}}
\nc{\scrB}{\mathscr{B}}
\nc{\scrC}{\mathscr{C}}
\nc{\scrD}{\mathscr{D}}
\nc{\scrF}{\mathscr{F}}
\nc{\scrK}{\mathscr{K}}
\nc{\A}{\D\forall}
\nc{\B}{\mathrm{B}}
\nc{\cB}{\mathcal{B}}
\nc{\cZ}{\mathcal{Z}}
\nc{\bB}{\mathbf{B}}
\nc{\BS}{\mathbf{B}(\mathcal{S})}
\nc{\BF}{\mathbf{B}(\mathcal{F})}
\nc{\BU}{\mathbf{B}(\mathcal{U})}
\nc{\cSp}{\mathcal{S}^+}
\nc{\cFp}{\mathcal{F}^+}
\nc{\cUp}{\mathcal{U}^+}
\nc{\BG}{\B_\Ga}
\nc{\BL}{\B_\Lambda}
\nc{\BT}{\B_\Tau}
\nc{\BTstar}{\B_{\Tau^*}}
\nc{\BO}{\B_\Om}
\nc{\DO}{\cD_\Om}
\nc{\KO}{\cK_\Om}
\nc{\CG}{C_\Ga}
\nc{\CL}{C_\Lambda}
\nc{\CT}{C_\Tau}
\nc{\CTstar}{C_{\Tau^*}}
\nc{\CO}{C_\Om}
\nc{\COgp}{C_{\Om^{\grbl}}}
\nc{\CLgp}{C_{\Lambda^{\grbl}}}
\nc{\BOgp}{\B_{\Om}^{\grbl}}
\nc{\BLgp}{\B_{\Lambda^{\grbl}}}
\nc{\sfC}{\mathsf{C}}
\nc{\sfD}{\mathsf{D}}
\nc{\bD}{\mathbf{D}}
\nc{\Tau}{\mathrm{T}}
\nc{\cA}{\mathcal{A}}
\nc{\cK}{\mathcal{K}}
\nc{\cD}{\mathcal{D}}
\nc{\cE}{\mathcal{E}}
\nc{\cF}{\mathcal{F}}
\nc{\cS}{\mathcal{S}}
\nc{\cT}{\mathcal{T}}
\nc{\cG}{\mathcal{G}}
\nc{\cY}{\mathcal{Y}}
\nc{\J}{\mathcal{J}}
\nc{\cL}{\mathcal{L}}
\nc{\cM}{\mathcal{M}}
\nc{\cN}{\mathcal{N}}
\nc{\cH}{\mathcal{H}}
\nc{\scF}{\mathscr{F}}
\nc{\scH}{\mathscr{H}}
\nc{\Op}{\mathrm{O}}
\nc{\rmA}{\mathrm{A}}
\nc{\rmF}{\mathrm{F}}
\nc{\rmB}{\mathrm{B}}
\nc{\rmD}{\mathrm{D}}
\nc{\rmP}{\mathrm{P}}
\nc{\cC}{\mathcal{C}}
\nc{\cP}{\mathcal{P}}
\nc{\bbQ}{\mathbb{Q}}
\nc{\bbR}{\mathbb{R}}
\nc{\cU}{\mathcal{U}}
\nc{\cX}{\mathcal{X}}
\nc{\cQ}{\mathcal{Q}}
\nc{\Un}{\bigcup}
\nc{\cV}{\mathcal{V}}
\nc{\cR}{\mathcal{R}}
\nc{\tcR}{\tilde{\mathcal{R}}}
\nc{\cW}{\mathcal{W}}
\nc{\Z}{{\mathbb Z}}
\nc{\Impl}{\Rightarrow}
\long\def\forget#1\forgotten{\marginpar{\textcolor{green}{Forgetting...}}}
\nc{\ft}{\mathfrak{t}}
\nc{\fb}{\mathfrak{b}}
\nc{\fc}{\mathfrak{c}}
\nc{\fd}{\mathfrak{d}}
\nc{\fg}{\mathfrak{g}}
\nc{\oo}{\infty}
\nc{\fr}{\mathfrak{r}}
\nc{\fk}{\mathfrak{k}}
\nc{\bidi}{\mathfrak{bidi}}
\nc{\fu}{\mathfrak{u}}
\nc{\fh}{\mathfrak{h}}
\nc{\fp}{\mathfrak{p}}
\nc{\fj}{\mathfrak{j}}
\nc{\fs}{\mathfrak{s}}
\nc{\w}{\omega}
\nc{\x}{\times}
\nc{\Iff}{\Leftrightarrow}
\newcommand\comp{^{\text{\tt c}}}
\nc{\nin}{\notin}
\nc{\cat}{\hat{\ }}
\nc{\sub}{\subseteq}
\nc{\spst}{\supseteq}
\nc{\sm}{\setminus}
\nc{\as}{\subseteq^*}
\nc{\les}{\le^*}
\nc{\leinf}{\le^{\infty}}
\nc{\leS}{\le_S}
\nc{\leF}{\le_F}
\nc{\leU}{\le_U}
\nc{\rest}{\restriction}

\nc{\la}{\langle}
\nc{\ra}{\rangle}
\nc{\E}{\exists}
\nc{\dom}{\op{dom}}
\nc{\cov}{\op{cov}}
\nc{\add}{\op{add}}
\nc{\addmen}{\add(\Men{})}
\nc{\cof}{\op{cof}}
\nc{\cf}{\op{cf}}
\nc{\non}{\op{non}}
\nc{\unif}{\op{non}}
\nc{\COV}{\op{COV}}
\nc{\ADD}{\op{ADD}}
\nc{\COF}{\op{COF}}
\nc{\NON}{\op{NON}}
\nc{\impl}{\to}
\nc{\Lp}{\mathcal{L_\p}}
\nc{\Wlog}{without loss of generality}
\newtheorem{thm}{Theorem}[section]
\nc{\bthm}{\begin{thm}} \nc{\ethm}{\end{thm}}
\newtheorem{need}[thm]{Need}
\nc{\bneed}{\begin{need}\color{dg}} \nc{\eneed}{\end{need}}
\newtheorem{prop}[thm]{Proposition}
\nc{\bprp}{\begin{prop}} \nc{\eprp}{\end{prop}}
\newtheorem{fact}[thm]{Fact}
\nc{\bfct}{\begin{fact}} \nc{\efct}{\end{fact}}
\newtheorem{prob}[thm]{Problem}
\nc{\bprb}{\begin{prob}} \nc{\eprb}{\end{prob}}
\newtheorem{lem}[thm]{Lemma}
\nc{\blem}{\begin{lem}} \nc{\elem}{\end{lem}}
\newtheorem{app}[thm]{Application}
\nc{\bapp}{\begin{app}} \nc{\eapp}{\end{app}}
\newtheorem{claim}[thm]{Claim}
\nc{\bclm}{\begin{claim}} \nc{\eclm}{\end{claim}}
\newtheorem{cor}[thm]{Corollary}
\nc{\bcor}{\begin{cor}} \nc{\ecor}{\end{cor}}
\newtheorem{conj}[thm]{Conjecture}
\nc{\bcnj}{\begin{conj}} \nc{\ecnj}{\end{conj}}
\theoremstyle{definition}
\newtheorem{defn}[thm]{Definition}
\nc{\bdfn}{\begin{defn}} \nc{\edfn}{\end{defn}}
\newtheorem{obs}[thm]{Observation}
\nc{\bobs}{\begin{obs}} \nc{\eobs}{\end{obs}}
\theoremstyle{remark}
\newtheorem{rem}[thm]{Remark}
\nc{\brem}{\begin{rem}} \nc{\erem}{\end{rem}}
\newtheorem{cnv}[thm]{Convention}
\nc{\bcnv}{\begin{cnv}} \nc{\ecnv}{\end{cnv}}
\newtheorem{exam}[thm]{Example}
\nc{\bexm}{\begin{exam}} \nc{\eexm}{\end{exam}}
\nc{\bpf}{\begin{proof}} \nc{\epf}{\end{proof}
}
\nc{\be}{\begin{enumerate}}
\nc{\ee}{\end{enumerate}}
\nc{\bi}{\begin{itemize}}
\nc{\bimy}{\my{\begin{itemize}}
\nc{\eimy}{\end{itemize}}}
\nc{\itm}{\item}
\nc{\ei}{\end{itemize}}
\nc{\Subsection}[1]{\goodbreak\subsection*{#1}}
\nc{\sone}{\mathsf{S}_1}
\nc{\sfin}{\mathsf{S}_\mathrm{fin}}
\nc{\ufin}{\mathsf{U}_\mathrm{fin}}
\nc{\Split}{\mathsf{Split}}

\nc{\gone}{\mathsf{G}_1}    
\nc{\tgfin}{{\mathsf{G}}^*_\mathrm{fin}}
\nc{\gfin}{\mathsf{G}_\mathrm{fin}}

\nc{\men}[1]{\sfin(\Op(#1),\Op(#1))}
\nc{\sch}{\ufin(\cO,\Omega)}
\nc{\rothb}{\text{Rothberger}}
\nc{\pmen}{\sfin(\Omega,\Omega)}
\nc{\Rothb}{\sone(\Op,\Op)}

\nc{\prothb}{\sone(\Omega,\Omega)}
\nc{\tU}{{\tilde{U}}}
\nc{\tF}{{\tilde{F}}}
\nc{\tY}{{\tilde{Y}}}
\nc{\tX}{{\tsup[1]{X}}}
\nc{\dtX}{{\tsup[2]{X}}}
\nc{\dt}[1]{{\tsup[2]{#1}}}
\nc{\td}{{\tilde{d}}}
\nc{\tz}{{\tilde{z}}}
\nc{\cfd}{\cf(\fd)}
\nc{\Gaprim}{\Gamma'}
\nc{\Gammaprim}{\Gaprim}
\nc{\Omprim}{\Omega'}
\nc{\Omegaprim}{\Omprim}
\nc{\Omctblprim}{\Omctbl'}

\nc{\msep}{\sfin(\cD,\cD)}
\nc{\rsep}{\sone(\cD,\cD)}
\nc{\cft}{\sfin(\Omega_{\mathbf{0}},\Omega_{\mathbf{0}})}
\nc{\scft}{\sone(\Omega_{\mathbf{0}},\Omega_{\mathbf{0}})}

\nc{\Umen}{U\text{-Menger}}
\nc{\hur}{\ufin(\cO,\Gamma)}
\nc{\tUmen}{\tU\text{-Menger}}

\nc{\Men}{\text{Menger}}
\nc{\Sch}{\text{Scheepers}}

\nc{\aspst}{\prescript{*}{}{\spst}\ }
\nc{\eqs}{=^*}
\nc{\ctblOm}{\Omega_{\mathrm{ctbl}}}
\nc{\GNga}{{\smallbinom{\Om}{\Ga}}}
\nc{\ctblga}{\smallbinom{\ctblOm}{\Ga}}

\nc{\nadd}{\cN_{\mathrm{add}}}
\nc{\ball}{\mathrm{B}}
\nc{\cOunif}{\cO^{\textrm{unif}}}

\nc{\sep}{
\vspace{2cm}
\noindent
\begin{minipage}{\textwidth}
	\textcolor{red}{\rule{\textwidth}{1pt}}
\end{minipage}
}
\nc{\FS}{\op{FS}}
\nc{\sums}{\op{SS}}
\nc{\SG}{\op{SG}}
\nc{\tSG}{\SG_\odot}
\nc{\G}{\op{G}}
\nc{\FSG}{\op{FSG}}
\nc{\FP}{\op{FP}}
\nc{\nonNadd}{\non(\nadd)}
\nc{\borga}{\Ga_\mathrm{Bor}}
\nc{\borw}{\Omega_\mathrm{Bor}}
\nc{\pick}{x}
\nc{\gen}{y}
\nc{\nullzind}{\sone(\{\Op_n^{\mathsf{unif}}\}_{n\in\bbN},\Ga)}
\nc{\nullzindf}[1]{\sone(\{\Op_{#1}^{\mathsf{unif}}\}_{n\in\bbN},\Ga)}

\definecolor{dg}{RGB}{42,101,24}

\nc{\myb}[1]{{\color{blue}{#1}}}
\nc{\mydg}[1]{{\color{dg}{#1}}}
\DeclareMathOperator{\eexists}{\exists}
\DeclareMathOperator{\fforall}{\forall}

\nc{\Exists}[1]{\bigl(\eexists #1\bigr)}
\nc{\Forall}[1]{\bigl(\fforall #1\bigr)}
\nc{\End}[1]{\bigl(#1\bigr)}
\nc{\dmo}[2]{\DeclareMathOperator{#1}{#2}}
\dmo{\Asc}{Asc}
\nc{\plusmin}{\wedge}

\nc{\cBsub}{{\cB^{\mbox{\tiny $\sub$}}}}

\nc{\Alice}{{\textsc{Alice}{}}}
\nc{\Bob}{{\textsc{Bob}}}

\nc{\bfO}{\mathbf{0}}

\nc{\proba}[1]{F}

\nc{\opn}[1]{\Op}
\nc{\sfinoo}[1]{\sfin(\Op,\Op)}
\nc{\om}[1]{\Om}

\nc{\restrict}{\upharpoonright}
\nc{\sgom}{\mathsf{S}_1(\Gamma,\Omega)}

\newcommand{\sbgg}{{\mathsf{S}_1(\Ga_\mathrm{Bor},\Ga_\mathrm{Bor})}}
\newcommand{\sbomom}{{\mathsf{S}_1(\Om_\mathrm{Bor},\Om_\mathrm{Bor})}}
\newcommand{\sbomfomf}{{\mathsf{S}_1(\Om_\mathrm{Bor}^\mathrm{fat},\Om_\mathrm{Bor}^\mathrm{fat})}}

\nc{\NNbar}{\overline{\w}^{\,\w}}

\usepackage[all]{xy}

\newlength{\spacebox}
\usetikzlibrary{arrows}
\DeclareMathOperator{\Fn}{Fn}
\DeclareMathOperator{\pr}{prod}

\author[M.~Pawlikowski]{Micha\l{} Pawlikowski}\address{Micha\l{} Pawlikowski, Faculty of Technical Physics, Information
Technology and Applied Mathematics,
Lodz University of Technology, Aleje Politechniki 8,
93-590 \L\'od\'z, Poland}
\email{michal-pawlikowski4@wp.pl}

\author[P. Szewczak]{Piotr Szewczak}
\address{Institute of Mathematics, University of Warsaw, Banacha 2,
02–097 Warsaw, Poland}
\email{p.szewczak@wp.pl}
\urladdr{http://piotrszewczak.pl}

\author[L.~Zdomskyy]{Lyubomyr Zdomskyy}
\address{Lyubomyr Zdomskyy, Institut f\"ur Diskrete Mathematik und Geometrie, Technische Universit\"at Wien, Wiedner Hauptstraße 8--10/104, 1040 Wien, Austria.}
\email{lzdomsky@gmail.com}
\urladdr{https://dmg.tuwien.ac.at/zdomskyy/}

\makeatletter
\@namedef{subjclassname@2020}{\textup{2020} Mathematics Subject Classification}
\makeatother

\subjclass[2020]{Primary: 05D10, 54D20, 54C35; Secondary: 16W22}
		
\keywords{scales, products, Scheepers diagram, Menger, Hurewicz, Rothberger}

\title{Scales, products and the second row of the Scheepers diagram}

\begin{document}
\begin{abstract}
We consider products of sets of reals with a combinatorial structure based on scales parameterized by filters.
Sets of this kind have been extensively investigated
in products of spaces with combinatorial covering properties such as the Hurewicz, Scheepers, Menger, and Rothberger properties.
We will complete this picture by focusing on properties from the second row of the Scheepers diagram.
We solve in the affirmative a problem due to Bartoszy{\'n}ski--Tsaban, whether all finite powers of a $\fb$-scale set satisfy the strong covering property $\sgw$.
We also show that, in the Miller model, the product of finitely many $\fd$-concentrated sets satisfies $\sgw$.
We further provide product counterexamples demonstrating limitations of the methods used.
\end{abstract}

\maketitle

\section{Introduction}
By \emph{space} we mean an infinite Tychonoff topological space.
For a class $\cA$ of covers of spaces and a space $X$, we denote by $\cA(X)$ the family of all covers of $X$ from the class $\cA$.
Let $\cA$ and $\cB$ be classes of covers of spaces.
We define the following combinatorial covering properties that a space $X$ can possess.
\begin{labeling}{$\ufin(\cA,\cB)$:}
	\item [$\sone(\cA,\cB)$:] for each sequence of covers $\eseq{\cU}\in\cA(X)$, there are sets $U_0\in\cU_0, U_1\in\cU_1,\dotsc$ such that $\set{U_n}{n\in\w}\in \cB(X)$; 
	\item [$\sfin(\cA,\cB)$:] for each sequence of covers $\eseq{\cU}\in\cA(X)$, there are finite sets $\cF_0\sub\cU_0,\cF_1\sub\cU_1,\dotsc$ such that $\Un_{n \in \w}\cF_n\in \cB(X)$;
	\item [$\ufin(\cA,\cB)$:] for each sequence of covers $\eseq{\cU}\in\cA$, none with a finite subcover, there are finite sets $\cF_0\sub\cU_0,\cF_1\sub\cU_1,\dotsc$ such that $\set{\Un\cF_n}{n\in\w}\in \cB$.	
\end{labeling}

Let $\cU$ be a cover of a space $X$.
The family $\cU$ is an \emph{$\w$-cover} of $X$ if $X \notin \cU$ and for each finite set $F \sub X$ there is a set $U \in \cU$ such that $F \sub U$.
The family $\cU$ is a \emph{$\gamma$-cover} of $X$ if $X \notin \cU$, it is infinite, and for each element $x \in X$ the set 
$\set{U \in \cU}{x \notin U}$ is finite.
By $\Opn, \Omega,$ and $\Ga$ we denote, respectively, the classes of all open covers, all open $\omega$-covers, and all open $\gamma$-covers.
For $\cA,\cB\in\{\Opn,\Om,\Ga\}$, the following diagram, called the Scheepers diagram~\cite{coc1}, shows the relations between combinatorial covering properties in the realm of Lindel\"of spaces.
  \begin{figure}[h]{
  \adjustbox{scale = 0.59, center}{
            \begin{tikzcd}[ampersand replacement=\&,column sep=.5cm]
	      \&\ \   \text{Hurewicz}\&\&\&\text{Scheepers}\& \text{Menger}\\[-1.5em]
	      {} \& \ \ \ufin(\Opn, \Gamma) \arrow{rr}\&\&\sfin(\Gamma,\Omega) \arrow{r}\arrow[leftarrow]{dd} \&\ufin(\Opn, \Omega)\arrow{r}\& \sfin(\Opn,\Opn)  \\
	     \ \ \sone(\Gamma,\Gamma)\arrow{rr}\arrow{ru}\&\&
	     \sone(\Gamma,\Omega)\arrow[crossing over] {rr}\arrow{ru}\& \&\sone(\Gamma,\Opn)\arrow{ru}\&\\
            \& \&\& \sfin(\Omega,\Omega) \&\&\\
	    \sone(\Omega,\Gamma)\arrow{rr}\arrow{uu}\&\ \ \ \ \ \ \& \sone(\Omega,\Omega)\arrow{rr}\arrow{uu}\arrow{ru}\&\& \sone(\Opn,\Opn)\arrow{uu}\ \&\\ [-1.5em]
	    {\gamma\text{-property}}\&\&  \&\&{\text{Rothberger}}\&
        \end{tikzcd}
  }
    }
\end{figure}

The mainstream of research on properties in the Scheepers diagram concerns \emph{sets of reals}, that is, spaces homeomorphic to subspaces of the Cantor space $\PN$.
By a \emph{set} with a topological property we mean a set of reals with this property.
By $\roth$ and $\Fin$ we denote, respectively, the families of all infinite and all finite subsets of $\w$.
We identify each $a \in \roth$ with its increasing enumeration, viewed as an element of the Baire space $\w^\w$.
The topologies inherited from $\Pof(\w)$ and $\NN$ coincide on $\roth$.
For functions $a, b \in \roth$, we write $a\les b$ if the set $\sset{n}{a(n)\leq b(n)}$ is cofinite.
A set $A\sub\roth$ is \emph{bounded} if there is a function $b\in\roth$ such that $a\les b$ for all $a\in A$.
A subset of $\roth$ is \emph{unbounded} if it is not bounded.
Let $\fb$ be the minimal cardinality of an unbounded subset of $\roth$.
The value of $\fb$ does not change if we consider its counterpart in $\NN$ instead of $\roth$.

A set $\set{x_\alpha}{\alpha<\fb} \sub \roth$ is a \emph{$\fb$-scale} if it is unbounded and $x_\alpha \le^\ast x_\beta$ for all $\alpha<\beta<\fb$.
A set is a \emph{$\fb$-scale set} if it is equal to $X\cup\Fin$ for some $\fb$-scale $X$.
Bartoszy\'{n}ski--Shelah~\cite{BS01} proved that any $\fb$-scale set (which exists in $\ZFC$) satisfies the Hurewicz property $\ufin(\Opn, \Gamma)$.
It was the first uniform $\ZFC$ counterexample to the Hurewicz conjecture~\cite{H1926}, which asserts that a set of reals has the Hurewicz property $\ufin(\Opn,\Gamma)$ if and only if it is a countable union of its compact subsets.
Sets of this kind play an important role when considering products with combinatorial covering properties.

For a topological property $\bfP$, a space $X$ is \emph{productively $\bfP$ in a given class of spaces} if, for every space $Y$ in that class having property $\bfP$, the product $X\x Y$ has property $\bfP$.
By a result of Miller, Tsaban, and Zdomskyy, each $\fb$-scale set is productively Hurewicz in the class of sets of reals~\cite[Theorem 6.5]{SPMProd}.
In particular, all finite powers of a $\fb$-scale set are Hurewicz.
A subset $D$ of $\roth$ is \emph{dominating} if for any function $a\in\roth$, there is a function $b\in D$ such that $a\les b$.
Let $\fd$ be the minimal cardinality of a dominating subset of $\roth$.
If $X\sub\roth$ is a $\fb$-scale that is dominating (such a set exists if and only if $\fb=\fd$~\cite{TsabanTheBook}), then the $\fb$-scale set $X\cup\Fin$ is productively Scheepers $\ufin(\Opn,\Omega)$ and productively Menger $\sfin(\Opn,\Opn)$ in the class of sets of reals~\cite[Theorem 6.2]{SPMProd}.
Let $\cov(\cM)$ be the minimal cardinality of a family of meager subsets of $\roth$ whose union is equal to $\roth$.
Assuming that $\cov(\cM)=\fb$, all finite powers of a $\fb$-scale set satisfy the Rothberger property $\sone(\Opn,\Opn)$~\cite[Corollary 6.9 (2)]{SPMProd} and \cite[Theorem 19]{CoC7}.
The origin of our investigations lies in the following question due to Bartoszy{\'n}ski--Tsaban.

\bprb[{Bartoszy{\'n}ski--Tsaban~\cite[Problem~15(2)]{bats}}]\label{prb:bats}
Do all finite powers of a $\fb$-scale set satisfy $\sgom$?
\eprb

In the original formulation, Bartoszy{\'{n}}ski--Tsaban asked about the property $\sgo$ instead of $\sgw$.
However, by Proposition~\ref{prp:sgosgw} below, these two versions of the problem are equivalent.
We solve this problem in the affirmative and study it in the context of the following generalizations of $\fb$-scale sets.

A \emph{semifilter} is a nonempty set $F \sub \roth$ such that for all $a \in F$ and $b \in \roth$ with $\card{a\sm b}<\w$, we have $b \in F$.
Examples of semifilters are the Frechet filter $\cFin$ of all cofinite subsets of $\w$, the full semifilter $\roth$, and every nonprincipal ultrafilter on $\w$.
For a semifilter $F$ and functions $a, b \in \roth$, we write $a \le_F b$ if the set $\sset{n}{a(n)\leq b(n)}$ is in $F$.
The relation $\les$ is the same as $\leq_{\cFin}$.
A set $X \sub \roth$ is \emph{$\le_F$-bounded} if 
there is $b \in \roth$ such that $a\le_F b$ for all $a\in X$.
A subset of $\roth$ is \emph{$\le_F$-unbounded} if it is not $\le_F$-bounded.
Let $\fb(F)$ be the smallest cardinality of a $\le_F$-unbounded subset of $\roth$.
Using this terminology, we have $\bof(\cFin)=\fb$ and $\bof(\roth)=\fd$.
Let $F$ be a semifilter.
A set $\sset{x_\alpha}{\alpha<\bof(F)}\sub\roth$ is a \emph{$\bof(F)$-scale}~\cite[Definition 2.8]{sfh} if it is $\le_F$-unbounded and $x_\alpha\le_F x_\beta$ for all $\alpha<\beta<\bof(F)$.
Then a $\bof(\cFin)$-scale is just a $\fb$-scale.
If $X\sub\roth$ is a $\bof(F)$-scale for some semifilter $F$, then we call the set $X\cup\Fin$ a \emph{$\bof(F)$-scale set}.

$\bof(F)$-scale sets and their variations  serve as important examples in the study of  products of sets with the Menger, Rothberger and other selection properties.
Every $\bof(F)$-scale is
an $F$-scale (the latter is recalled in Section~\ref{sec:applications}), and a good part of our results holds already for $F$-scales.
More precisely, for each filter $F$, all finite powers of a $\bof(F)$-scale set are Menger~\cite[Theorem 4.1]{sfh}.
Moreover, if $U$ is an ultrafilter with $\fb(U)\leq\cov(\cM)$, then all finite powers of a $\bof(U)$-scale set are Rothberger~\cite[Lemma 2.21]{stz}.
By a result of Szewczak, Tsaban, and Zdomskyy~\cite[Theorem 2.5, Lemma 2.10]{stz}, under the assumptions that $\fd\leq\fr$ and that $\fd$ is regular, there exist a $\bof(U)$-scale set and a 
$\bof(\tilde{U})$-scale set for some ultrafilters $U$ and $\tilde{U}$, whose product space is not Menger.
In particular, this shows that, under these assumptions, there are two Menger sets that are Menger in all finite powers (equivalently, they satisfy $\sfinww$), but whose product is not Menger.
The results in~\cite{stz} cited above are formulated in terms of $U$-scales (see Section~\ref{sec:applications} for more details). 
For~\cite[Lemma 2.21]{stz}, this causes no difficulty, since every $\bof(U)$-scale is a $U$-scale. 
Although the above-mentioned results~\cite[Theorem 2.5, Lemma 2.10]{stz} are stated for a $U$-scale and a $\tilde U$-scale, respectively, a small modification of their construction yields a $\bof(U)$-scale and a $\bof(\tilde U)$-scale.
The above results show that $\bof(U)$-scale sets, where $U$ is an ultrafilter, serve as an important source of counterexamples to productivity of combinatorial covering properties, in particular for the Menger property and for local properties in function spaces~\cite[Proposition 3.1]{stz}.

Let $\borga$ and $\borw$ be the classes of all \emph{countable Borel $\gamma$-covers} and \emph{countable Borel $\w$-covers} of spaces, respectively.
The common combinatorial feature behind the scale constructions used in our positive results is isolated in Definition~2.1. For an infinite cardinal $\kappa$, a set $X\sub\roth$ is called \emph{$\kappa$-fin-unbounded} if $\card{X}\geq\kappa$ and, for every $d\in\roth$, there is a set $S\sub X$ with $\card{S}<\kappa$ such that every finite set $H\sub X\sm S$ has the property that $\bigcup H$ omits infinitely many intervals $[d(n),d(n+1))$. We call $X\cup\Fin$ a \emph{$\kappa$-fin-unbounded scale set}. This condition is designed to control finite subsets of a scale simultaneously and is therefore particularly well suited to product arguments. In Section~\ref{sec:applications} we connect it with the filter-parametrized scales discussed above: Proposition~\ref{lem:F-scale_is_b(f)_unbdd} shows that, for every filter $F$, every $F$-scale is a $\fb(F)$-fin-unbounded set.

Our first main result gives a uniform product theorem for these sets and provides a $\ZFC$ solution to Problem~\ref{prb:bats}.
{
\renewcommand{\thethm}{2.2}
\bthm
Let $\kappa$ be an infinite regular cardinal number with $\kappa\leq\fd$, let $X$ be a $\kappa$-fin-unbounded set, let $k$ be a natural number, and let $Y$ be a set satisfying $\sbgg$. Then the product space $(X\cup\Fin)^k\x Y$ satisfies $\sgom$.
\ethm
}

The proofs of this result and its applications use the Michael topology on the relevant scale sets as an auxiliary device. Isolating the infinite points in this topology makes it possible to streamline the combinatorial arguments for products; the statements displayed here are the resulting versions for the standard topology. Together with the scale results in Section~\ref{sec:applications}, Theorem~2.2 yields the product results for $\bof(F)$-scale sets mentioned above.

We also obtain a corresponding theorem involving the classical notion of a $\gamma$-set, introduced by Gerlits and Nagy~\cite{GN1982}. A set of reals $Y$ is a \emph{$\gamma$-set} if every open $\omega$-cover of $Y$ contains a $\gamma$-subcover; equivalently, the set $Y$ satisfies $\sone(\Omega,\Gamma)$.
{
\renewcommand{\thethm}{4.1}
\bthm
Let $\kappa$ be an infinite regular cardinal with $\kappa\leq\cov(\cM)$, let $X$ be a $\kappa$-fin-unbounded set, let $k$ be a natural number, and let $Y\sub\PN$ be a $\gamma$-set. Then the product space $(X\cup\Fin)^k\x Y$ satisfies $\sone(\Omega,\Omega)$.
\ethm
}
This result extends earlier product theorems for scale sets and illustrates that the same combinatorial framework applies to several properties from the second and third rows of the Scheepers diagram.

The positive results have natural limitations. In Section~\ref{sec:counterexamples} we show that in Theorem~\ref{thm:kappa_fin_unb_x_sbgg_is_sgom} the assumption $\sbgg$ cannot in general be replaced by $\sbomom$. The first counterexample is obtained by combining ideas of Szewczak and Wi\'{s}niewski~\cite{luzincomb} with those of Bartoszy\'{n}ski, Shelah, and Tsaban~\cite{BST2003}.
{
\renewcommand{\thethm}{5.1}
\bthm
Assume that $\cov(\cM)=\fc$ and $\fc$ is regular. Then there exist ultrafilters $U$ and $\tilde U$, a $U$-scale $X$, and a $\tilde U$-scale $Y$, such that $X\cup\Fin$ and $Y\cup\Fin$ satisfy $\sbomom$, but the product space $(X\cup\Fin)\x(Y\cup\Fin)$ is not Menger. In fact, even $X$ and $Y$ satisfy $\sbomom$, and $X\x Y$ is not Menger.
\ethm
}
A second construction, using Cohen forcing, shows that this obstruction already occurs for a classical $\fb$-scale.
{
\renewcommand{\thethm}{5.8}
\bthm
It is consistent with $\CH$ that there exist a $\fb$-scale set and a set $Y$ satisfying $\sbomom$ whose product space is not Menger.
\ethm
}
These two theorems delimit the extent to which the positive product results above can be strengthened.

Finally, we apply the same framework to concentrated sets in the Miller model. A set of reals $X$ is \emph{$\fd$-concentrated on a countable subset $D$ of $X$} if $\card{X}\geq\fd$ and $\card{X\sm U}<\fd$ for every open set $U$ containing $D$. In $\ZFC$ there exists a set that is $\fd$-concentrated on $\Fin$, and each such set satisfies $\sgo$~\cite[Corollary 1.14]{TsabanTheBook}. Miller, Tsaban and Zdomskyy proved that in the Miller model the product of finitely many sets that are $\fd$-concentrated on $\Fin$ is Menger~\cite[Theorem~4.7]{SPMProd}. For $X\sub\roth$, being \emph{$\fd$-unbounded} means that $\card{X}\geq\fd$ and $\card{\sset{x\in X}{x\leq c}}<\fd$ for every $c\in\roth$; equivalently, $X\cup\Fin$ is $\fd$-concentrated on $\Fin$~\cite[Lemma 2.3]{ST}. Under the Near Coherence of Filters principle ($\NCF$), these sets have a particularly useful scale structure.
{
\renewcommand{\thethm}{6.2}
\bthm
Assume that $\NCF$ holds. Let $U\sub\roth$ be an ultrafilter. Each $\fd$-unbounded set in $\roth$ is a $U$-scale.
\ethm
}

Since $\NCF$ holds in the Miller model, Theorem~\ref{thm:fd_unbdd_in_Miller} together with the preceding product results gives the following strengthening of the Menger-product theorem.
{
\renewcommand{\thethm}{6.3}
\bcor
In the Miller model, for every positive natural number $k$ and $\fd$-concentrated sets $Y_i\sub\PN$ on  
$\fin$ for $i<k$, the product $\prod_{i
< k} Y_i$ satisfies $\sgw$.
\ecor
}
We are grateful to the referee for the careful reading of the manuscript and for the helpful comments and suggestions, which led us to improve substantially the presentation of the paper, especially the Introduction.

\section{Finite powers of $\kappa$-fin-unbounded scale sets}

For natural numbers $k,l$ with $k<l$, let $[k,l):=\sset{i}{k\leq i <l}$.
We say that a set $x$ \emph{omits an interval of $a \in \roth$} if $x \cap [a(n), a(n+1)) = \emptyset$ for some natural number $n$.

\bdfn
Let $\kappa$ be an infinite cardinal number.
A set $X \sub \roth$ is a \emph{$\kappa$-fin-unbounded set} if $\card{X} \geq \kappa$ and for each $d \in \roth$
there is a set $S \sub X$ with $\card{S}<\kappa$ such that for every finite set $F \sub X \sm S$ the set $\bigcup F$ omits infinitely many intervals of $d$. 
A set $X\cup\Fin\sub\PN$ is an \emph{$\kappa$-fin-unbounded scale set} if $X$ is a $\kappa$-fin-unbounded set.
\edfn

The \emph{Michael topology}~\cite{M1963} on $\PN$ is the topology in which open neighborhoods of points in $\Fin$ are the same as in the standard topology on $\PN$ and all points of $\roth$ are isolated. 
For $X \sub \PN$, we denote by $X_\mathrm{M}$ the set $X$ endowed with the Michael subspace topology.

The main goal of this section is to prove the following result.

\bthm\label{thm:kappa_fin_unb_x_sbgg_is_sgom}
Let $\kappa$ be an infinite regular cardinal number with $\kappa \leq \fd$, let $X$ be a $\kappa$-fin-unbounded set, let $k$ be a natural number, and let $Y$ be a set satisfying $\sbgg$.

\begin{enumerate}
\item The space
$\michk{X \cup \Fin}{k} \x Y$ satisfies $\sgom$.
\item The space $(X \cup \Fin)^{k} \x Y$ satisfies $\sgom$.
\end{enumerate}
\ethm

In Section~\ref{sec:applications} we give equivalent conditions for being a $\kappa$-fin-unbounded {set} and also show that this notion captures many subsets of $\roth$ with a combinatorial structure, e.g., finite unions of $\fb$-scales or finite unions of $U$-scales for some ultrafilter $U$.
In Section~\ref{sec:applications} we also give applications of Theorem~\ref{thm:kappa_fin_unb_x_sbgg_is_sgom}.
In particular, we obtain several results on products of subspaces of $\PN$ 
with the standard topology.

In order to prove Theorem~\ref{thm:kappa_fin_unb_x_sbgg_is_sgom}, we need the following notion and auxiliary results.
For a nonempty finite set $F\sub\roth$, let $\min F\in\roth$ be a function such that $(\min F)(n) \coloneqq \min \set{x(n)}{x \in F}$ for all $n$.
For functions $a,b\in\roth$ we write $a\leinf b$ if the set $\sset{n}{a(n)\leq b(n)}$ is infinite.

Let $k$ be a natural number and let $Y$ be a space.
The \emph{canonical topology on} $\PN^k \x Y$ is the product topology obtained by equipping 
$\PN$ with its standard topology.
For a set $U \sub \PN^k \x Y$, let $\Intc(U)$ denote the interior of $U$ with respect to the canonical topology on $\PN^k \x Y$.
Let $X\sub\roth$.
Whenever $\cU$ is a family of open subsets of $\michk{X\cup\Fin}{k}\x Y$, we choose a family $\cU'$ of open subsets of $\PN_\mathrm{M}^k\x Y$ such that
\[
\cU=\sset{U\cap((X\cup\Fin)^k\x Y)}{U\in\cU'}.
\]
In our considerations, we will not distinguish between these two families and will usually regard $\cU$ as a family of open sets in the ambient space $\PN_\mathrm{M}^k\x Y$.
We use the analogous convention when the spaces are endowed with their standard topologies and for other kinds of covers.

\bfct\label{fct:interior_of_gamma_cover_in_MTOP}
Let $k$ be a natural number, let $X\sub\roth$, and let $Y$ be a space.
Then the following assertions hold.
\be
\itm
For each open set $U\sub\PN_\mathrm{M}^k\x Y$, we have
\[
U\cap(\Fin^k\x Y)\sub\Intc(U).
\]

\itm
Let $\sset{U_n}{n\in\bbN}\in\Gamma(\michk{X\cup\Fin}{k}\x Y)$.
Then, for each $(x,y)\in\Fin^k\x Y$, we have
\[
(x,y)\in\Intc(U_n)
\]
for all but finitely many $n$.

\itm 
Let $\eseq{U}$ be a sequence of open sets in 
$\PN_{\mathrm{M}}^k \x Y$ 
such that for each $(x,y) \in \Fin^k \x Y$, we have $(x,y) \in U_n$ for all 
but finitely many $n$ and $\Intc(U_n) \neq \emptyset$ for all $n$.
Then there exists a family of pairwise distinct sets $\sset{V_n}{n \in \w}$ open in $\PN^k \x Y$ such that $V_n\sub U_n$ for all $n\in\omega$, 
and for each $(x,y) \in \Fin^k \x Y$, we have $(x,y) \in V_n$ for all but finitely
many $n$.
\ee
\efct

\bpf
(1)~Take any $(x,y)\in U\cap(\Fin^k\x Y)$.
Then there are open sets $U_0,U_1,\dotsc,U_{k-1}\sub\PNM$ and an open set $U_y\sub Y$ such that
\[
(x,y)\in\prod_{i<k}U_i\x U_y\sub U.
\]
Since $x\in\Fin^k$, there are open sets $U_0',U_1',\dotsc,U_{k-1}'\sub\PN$ such that
\[
x\in\prod_{i<k}U_i'
\quad\text{and}\quad
U_i'\sub U_i
\]
for all $i<k$.
Thus, the set $\prod_{i<k}U_i'$ is open in $\PN^k$ and
\[
(x,y)\in\prod_{i<k}U_i'\x U_y
\sub\prod_{i<k}U_i\x U_y
\sub U.
\]
Therefore, $(x,y)\in\Intc(U)$.

(2)~Fix $(x,y)\in\Fin^k\x Y$.
Then $(x,y)\in U_n$ for all but finitely many $n$.
By~(1), we have $(x,y)\in\Intc(U_n)$ for all but finitely many $n$.

(3)~By~(2), for each $(x,y) \in \Fin^k \x Y$, we have $(x,y) \in \Intc(U_n)$ for all 
but finitely many $n$.
Now by recursion we construct a sequence $\seq{x_n}{n \in \w}$ and sets 
$\eseq{V}$ open in $\PN^k \x Y$ such that $V_n \sub U_n$, $x_n \in V_n$, and $\sset{x_i}{i<n} \cap V_n = \emptyset$ for all $n$.
Define $V_0 \coloneqq \Intc(U_0)$. Since $V_0 \neq \emptyset$, there exists 
$x_0 \in V_0$. Now suppose we have defined $V_0, V_1, \dotsc, V_n$ and 
$x_0, x_1, \dotsc, x_n$.
Since $\Intc(U_{n+1}) \neq \emptyset$, the set $V_{n+1} \coloneqq \Intc(U_{n+1}) \sm \sset{x_i}{i\leq n}$ is non-empty and open.
Observe that since $x_i \notin V_{n+1}$, it follows that $V_i \neq V_{n+1}$ for all $i<n+1$. Pick $x_{n+1} \in V_{n+1}$.
\epf

Let $\overline{\w}:=\w\cup\{\infty\}$, where $\infty\notin\w$.
We extend the standard order $\leq$ on $\w$ to $\overline{\w}$ in a natural way, declaring that $n\leq \infty$ for all natural numbers $n$.
In $\overline{\w}$ we consider the discrete topology.
In the space $\NNbar$, define a relation $\les$ in the same way as in $\NN$.

By $\Omprim$ we denote the class of all families of open sets {such that
for each finite subset of the space, 
it is contained in infinitely many members of the family,} without the requirement that no member of the family contains the whole space.
More precisely, if $X$ is a subspace of a space $Y$, then $\Omprim(X)$ consists of all families $\cU$ of open subsets of $Y$ such that for each finite set $F\sub X$ there {exist infinitely many sets} $U\in\cU$ with $F\sub U$.
By $\Gaprim$ we denote the class of all infinite families of open sets satisfying the point-cofinite condition, without the requirement that no member of the family contains the whole space.
More precisely, if $X$ is a subspace of a space $Y$, then $\Gaprim(X)$ consists of all infinite families $\cU$ of open subsets of $Y$ such that, for each $x\in X$, the set $\sset{U\in\cU}{x\notin U}$ is finite.
The ambient space for $X$ will always be clear from the context.
We will use the same convention for analogous variants of $\Omega$ and $\Gamma$.

\blem[{\cite[Lemma~3.2]{unbddtower}}]
\label{lem:unbddtower}
Let $m$ be a natural number, let $Y$ be a set, and let $\sset{U_n}{n \in \omega}$ be a family of open sets in $\PN^m \x Y$ such that for each point $(x,y)\in \Fin^m\x Y$, we have $(x,y)\in U_n$ for all but finitely many $n$. 
There are Borel   functions $f,g \colon Y \to \NNbar$ and a Borel function $d \colon Y \to \NN$ such that the following assetions hold.
\be
\item We have $f(y)^{-1}(\infty)= g(y)^{-1}(\infty)$, and this preimage is finite.
\item The functions $f(y)\upharpoonright \big(\omega\setminus f(y)^{-1}(\infty)\big)$
and $g(y)\upharpoonright \big((\omega\setminus f(y)^{-1}(\infty)\big)$ are non-decreasing
and finite-to-one, and $f(y)(n)\leq g(y)(n)$ for all $n\in\omega$.
\item For all $x \in (\roth)^m, y \in Y$, and $n \in \w$ 
 \[
\text{if }x_i \cap [f(y)(n), g(y)(n)) = \emptyset\text{ for all }i < m,\text{ then }(x,y) \in U_n.
\]
 \item We have
 \[
 f(y)(d(y)(i))=i
 \] 
for all but finitely many $i$. Consequently, the restriction of $d(y)$ to the cofinite set of all $i\in\omega$ for which the equality above holds is nondecreasing.
\ee
\elem
\brem
Lemma~\ref{lem:unbddtower} is only the preliminary part of the result due to Szewczak and Włudecka~\cite[Lemma~3.2]{unbddtower}.
Moreover, in the original formulation of Lemma~\ref{lem:unbddtower} we request that $\sset{U_n}{n\in\bbN}\in\Gamma(\Fin^m\x Y)$ but the proof allows to weaken the assumption to the above form.
\erem
\blem\label{lem:Holes_in_bscale_times_sbgg}
Let $\kappa$ be an infinite cardinal number, let $X\sub\roth$ be a $\kappa$-fin-unbounded set, and let $Y$ be a set satisfying $\sbgg$.
Let $A\sub X$ be a set with $\card{A}<\fb$, let $k$ be a natural number, and let $\cU\in\Gammaprim(\michk{\Fin\cup A}{k}\x Y)$ be a family of open sets in ${\PN_{\mathrm{M}}^k} \x Y$.
Then there is a set $S \sub X$ with $\card{S}<\kappa$ such that
{each finite set $F \sub \michk{A \cup \Fin \cup X \sm S}{k}\x Y$ is contained in some member of $\cU$.}
\elem

\bpf
Let $\cU=\sset{U_n}{n\in\w}$.
Fix $I\sub k$.
Let
\[
B_{I} \coloneqq \Fin^{k\setminus I} \x (\Fin \cup A)^I\x Y.
\]
For a set $V\sub \PN^{k\sm I}\x \PN_{\mathrm{M}}^I\x Y$, let $\Int^I_\mathrm{c}(V)$ be the interior of $V$ in the product topology on $\PN^{k\sm I}\x \PN_\mathrm{M}^I\x Y$.
By Fact~\ref{fct:interior_of_gamma_cover_in_MTOP}(2) we have that each element of $B_I$ belongs to $\Int^I_\mathrm{c}(U_n)$ for all but finitely many $n$.
Therefore, $\Int^I_{\mathrm{c}}(U_n) \neq \emptyset$ for all but finitely many $n$.
Let $n_I$ be the smallest $n$ with the above property.
By Fact~\ref{fct:interior_of_gamma_cover_in_MTOP}~(3), there exist a 
family $\cV_I \coloneqq \sset{V^I_n}{n \in \w}$ of sets open in $\PN^{k \sm I} \x (\Fin \cup A)^I \x Y$ such that $V^I_n\sub \Int^I_{\mathrm{c}}(U_n)$ for all $n \ge n_I$, and each member of $B_I$ belongs to $V^I_n$ for all but finitely many $n$.
Let 
\[
d_{I}, f_{I}, g_{I} \colon 
(\Fin \cup A)^I_M \x Y \to \NNbar
\]
be the Borel functions supplied by
Lemma~\ref{lem:unbddtower}, applied to the family {$\cV_I$} with $Z_I=\michk{\Fin \cup A}{I} \x Y$ in place of $Y$.
Since $\michk{\Fin \cup A}{I} \x Y$ satisfies $\sbgg$ for all $I\subseteq k$, there is a function $b \in \roth$ 
with $b(0) \neq 0$ such that
\[
 \set{d_I(\vec{a}, y),g_I (\vec{a}, y)}{(\vec{a}, y)\in \michk{\Fin \cup A}{I} \x Y, I\sub k}\les b.
\]
Let $c\in\roth$ be a function such that 
$c(0):=b(0)$ and $c(n+1):=b(c(n))+1$ for all natural numbers $n$.
Then for all $I\sub k$ and $(\vec{a}, y) \in (\Fin \cup A)^I_M \x Y$, we have
  \begin{equation}\label{eq:con4}
    c(m) \leq f_{I}(\vec{a}, y)(c(m+1))\leq g_{I}(\vec{a}, y)(c(m+1))<c(m+2)
  \end{equation}
for all but finitely many natural numbers $m$.
Indeed, the middle inequality follows from Lemma~\ref{lem:unbddtower}(2), the third inequality from 
$c(m+2)=b(c(m+1))+1$ and the fact that
\[
g_I(\vec{a}, y)\les b,
\]
while the first inequality can be established as follows:
\begin{eqnarray*} c(m)=f_I(\vec{a}, y)\Big(d_I(\vec{a}, y)\big(c(m)\big)\Big)\leq
f_I(\vec{a}, y)\Big(b\big(c(m)\big)\Big)
\le f_I(\vec{a}, y)\big(c(m+1)\big),
\end{eqnarray*}
where we used the monotonicity of
$f_I(\vec{a}, y)$ on a cofinite set 
as well as $d_I(\vec{a},y)\les b$.

Let $\tilde{c} \in \roth$ be a function such that $\tilde{c}(n):=c(2n)$ for all $n$.
Since $X$ is $\kappa$-fin-unbounded, there exists a set $S \sub X$ with $\card{S}<\kappa$ such that for each finite set $F \sub X \sm S$ the set $\Un F$ omits infinitely many intervals
of $\tilde{c}$.
Fix
\[
H\in [\michk{A \cup \Fin \cup X \sm S}{k}\x Y]^{<\omega}
\]
and find
$F_X\in [X\sm S]^{<\omega}$, $F_A\in [A\cup \Fin]^{<\omega}$
and $F_Y\in [Y]^{<\omega}$ such that 
\[
 H\sub\bigcup\sset{F_X^{k\setminus I}\times F_A^I\times F_Y}{I\sub k}.
\]
Then there is a natural number $m_0$ such that the inequalities in~\eqref{eq:con4}
hold for all
$I\sub k$, $(\vec{a}, y) \in F_A^I \x F_Y$ and $m \geq m_0$.
Now take any $m_1 \geq m_0$ such that $\Un F_X$ omits the interval 
\[
[\tilde{c}(m_1), \tilde{c}(m_1+1))=[c(2m_1), c(2m_1+2)).
\]

We have $H \sub U_{c(2m_1+1)}$. Indeed,
take any $h \in H$.
Then there are 
\[
I_h\sub k,\quad
\vec{x}=(x_j)_{j\in k\sm I_h}\in (X\sm S)^{k\sm I_h}\quad\text{and}\quad(\vec{a},y)\in F^{I_h}_A\times F_Y
\] 
such that $h\upharpoonright (k\sm I_h)=\vec{x}$, $h\upharpoonright I_h=\vec{a}$ and $y$ is the $Y$-coordinate of $h$.
We have 
\[
\bigcup_{j\in k\setminus I}x_j \cap [f_{I_h}(\vec{a}, y)(c(2m_1+1)), g_{I_h}(\vec{a}, y)(c(2m_1+1)))\sub\bigcup_{j\in k\setminus I}x_j \cap[c(2m_1),c(2m_1+2))=\emptyset,
\]
and thus 
$
h=(\vec{x}, \vec{a}, y) \in {V^{I_h}_{c(2m_1+1)}}.
$ 
It follows that
\[
H\sub\bigcup_{h\in H}V^{I_h}_{c(2m_1+1)} \sub U_{c(2m_1+1)},
\]
which completes the proof.
\epf

\bcor\label{cor:Holes_in_Fscale_times_Y}
Let $k$ be a natural number and $Y$ be a set satisfying $\sbgg$.
Let $\kappa$ be an infinite cardinal number, $X$ be a $\kappa$-fin-unbounded set, $H \sub X$ be a finite set and $\cU \in \Gammaprim(\michk{\Fin \cup H}{k} \x Y)$ be a family of open sets in $\PN_{\mathrm{M}}^k \x Y$.
Then there is a set $S \sub X$ 
with $\card{S}<\kappa$ such that 
{each finite set $F \sub \michk{H \cup \Fin \cup X \sm S}{k}\x Y$ is contained in some member of $\cU$.}
\ecor
\blem\label{lem:general_way_of_proof}
Let $X\sub\roth$, $Y$ be a space, $\kappa$ be an infinite regular cardinal number, $k$ be a natural number and $\cU_{n,m} \in \Opn(\michk{X \cup \Fin}{k} \x Y)$ such that no element of $\cU_{n,m}$ contains $\michk{X \cup \Fin}{k} \x Y$ for all natural numbers $n, m$.
Assume that for every set $A \sub X$ with $\card{A}<\kappa$ and every natural number $n$, there are sets $U_{n, m}(A) \in \cU_{n, m}$ for all $m\in\omega$ such that for each finite set $H \sub A$, there is a set $S_H \sub X$ with $\card{S_H}<\kappa$ {such that for each finite set $F \sub \michk{H \cup \Fin \cup X \sm S_H}{k} \x Y$, there exists $m \in \w$ such that 
$F \sub U_{n,m}(A)$.}
Then there exist $U_{n,m} \in \cU_{n,m}$ for all $n,m \in \w$ such that 
\[
\set{U_{n,m}}{n,m \in \w} \in \Omega(\michk{X \cup \Fin}{k} \x Y).
\]
\elem
\bpf
Let $A_{-1}:=\emptyset$. Using our assumption for 
$n=0$, and $A=A_{-1}$, we can find 
 $U_{0,m}(A_{-1}) \in \cU_{0,m}$ for all natural numbers $m$ and a set $S_0 \sub X$ with $\card{S_0}<\kappa$ 
{such that for each finite set $F \sub \michk{\Fin \cup X \sm S_0}{k} \x Y$, there exists $m \in \w$ such that 
$F \sub U_{0,m}(A_{-1})$.}
Set $A_0=S_0$ and fix a natural number $n$. 

{Assume that sets $U_{n,m}(A_{n-1}) \in \cU_{n,m}$ for all natural numbers $m$ and a set $A_n\sub X$ with $\card{A_n}<\kappa$ such that
\be 
\itm $A_{n-1} \sub A_{n}$;
\itm For each finite set $F \sub \michk{A_{n-1} \cup \Fin \cup X \sm A_n}{k} \x Y$, there exists $m \in \w$ such that $F \sub U_{{n-1},m}(A_{n-1})$;
\ee
have been defined.}
Then, by our assumption for $n+1$ and $A_n$ there are sets $U_{n+1,m}(A_n) \in \cU_{n+1,m}$ for all natural numbers $m$ such that for each finite set 
$H \sub A_n$ there is a set $S_H \sub X$ with $\card{S_H}<\kappa$ {such that
for each finite set $F \sub \michk{H \cup \Fin \cup X \sm S_H}{k} \x Y$,
there exists $m \in \w$ such that $F \sub U_{n+1,m}(A_n)$.}
Define 
\[
A_{n+1} \coloneqq A_{n} \cup \Un\sset{S_H}{H\sub A_{n}\text{ is finite}}.
\] 
Since $\kappa$ is regular, we conclude that $\card{A_{n+1}}<\kappa$.
{Therefore for each finite set \[F \sub \michk{A_{n} \cup \Fin \cup X \sm A_{n+1}}{k} \x Y,\] there exists $m \in \w$ such that 
$U_{n+1,m}(A_n)$.}

Fix finite sets $H_X \sub X \cup \Fin$ and $H_Y \sub Y$, and let $A \coloneqq \Un_{n \in \w} A_n$.
Since the sequence $\eseq{A}$ is increasing, there is a natural number $n$ such that $H_X\cap A=H_X\cap A_n$, i.e.,
\[
H_X \sub A_n \cup \Fin \cup X \sm A\sub A_n \cup \Fin \cup X \sm A_{n+1}.
\]
It follows that the set $H_X^k \x H_Y$ is covered by one of the sets $U_{n+1, m}(A_n)$. Set $U_{n,m}:=U_{n,m}(A_{n-1})$ for all $n,m\in\omega$.
\epf

\blem\label{lem:Quasi_omega_cover_of_AxSbgg}
Let $A$ be a space with $\card{A}<\fd$, let $k$ be a natural number, let $Y$ be a set satisfying $\sbgg$, and let 
$\cU_n=\sset{U_{n,m}}{m \in \omega} \in \Gammaprim(\michk{A \cup \Fin}{k} \x Y)$ {be a family of open sets in $\PN_{\mathrm{M}}^k \x Y$} for all $n$.
Then there are sets $U_n \in \cU_n$ for all $n$ 
such that for each finite set $F \sub A$ there is a family 
$\cU_F \sub \sset{U_n}{n \in \omega}$ such that $\cU_F \in \Gammaprim(\michk{F \cup \Fin}{k} \x Y)$.
\elem
\bpf
Fix a natural number $n$.
Let $V_{n,m} \coloneqq \bigcap_{i \geq m} U_{n,i}$ for all $m$. 
Then the family $\cV_n \coloneqq \sset{V_{n,m}}{m \in \omega}$ is an increasing {Borel cover} of $\michk{A \cup \Fin}{k} \x Y$.
Let $F \sub A$ be a finite set.
Then $\michk{F \cup \Fin}{k} \x Y$ satisfies $\sbgg$.
A function $\varphi_F \colon \michk{F \cup \Fin}{k} \x Y \to \omega^\omega$ given by 
\[
 \varphi_F(x,y)(n) \coloneqq \min \sset{m}{(x,y) \in V_{n,m}},
\]
for all $n$, is Borel. 
Since the space $\michk{F \cup \Fin}{k} \x Y$
is a countable union of topological copies of $Y$,
which is a set of reals,
by a result of Scheepers and Tsaban~\cite[Theorem 1]{ST2002}, the image of $\varphi_F$ is bounded in $\omega^\omega$.
Thus, 
there is a function $b_F \in \omega^\omega$ such that $\varphi_F[\michk{F \cup \Fin}{k} \x Y]\les b_F$.
Since $\card{A}<\fd$, there is a function $c \in \roth$ such that 
$\sset{b_F}{F\sub A\text{ is finite}} \leinf c$.
Choose $c$ in such a way that the sets $U_{n,c(n)}$ are pairwise different.
Fix a finite set $F \sub A$.
{For each $(x,y) \in \michk{F \cup \Fin}{k} \x Y$, $(x,y) \in V_{n, b_F(n)}$ 
for all but finitely many $n$.}
It follows that 
\[
\sset{U_{n, c(n)}}{b_F(n)<c(n), n\in\w} \in \Gammaprim(\michk{F \cup \Fin}{k} \x Y).\qedhere
\]
\epf

\brem 
Observe that the proof above uses no topological properties of $(A \cup \Fin)^k \x Y)$ other than the assumption that $Y$ satisfies $\sbgg$. Hence Lemma~\ref{lem:Quasi_omega_cover_of_AxSbgg} remains valid when $A \cup \Fin$ is endowed with any topology.
\erem

\bpf[{Proof of Theorem~\ref{thm:kappa_fin_unb_x_sbgg_is_sgom}}]
(1)
Let $\cU_{n,m} \in \Gamma(\michk{X \cup \Fin}{k} \x Y)$.
Fix a natural number $n$ and a set $A \sub X$ with $\card{A}<\kappa$.
By Lemma~\ref{lem:Quasi_omega_cover_of_AxSbgg}, there are sets $U_{n,m}(A) \in \cU_{n,m}$ for all $m$ such that 
for each finite set $H \sub A$ there is a family $\cU_H \sub \set{U_{n,m}(A)}{m \in \w}$ such that 
$\cU_H \in \Gammaprim({\michk{H \cup \Fin}{k}} \x Y)$ {and no element of $\cU_H$ contains $\PN^k \x Y$}.
Fix a finite set $H\sub A$.
By Corollary~\ref{cor:Holes_in_Fscale_times_Y}, there is a set $S_H \sub X$ with 
$\card{S_H}<\kappa$ such that 
{each finite set $F \sub \michk{\Fin \cup H \cup X \sm S_H}{k} \x Y$ is contained in some member of $\cU_H$.}
Thus each finite set $F \sub \michk{\Fin \cup H \cup X \sm S_H}{k} \x Y$ is contained in some member of $\set{U_{n,m}(A)}{m \in \w}$
for all finite sets $H \sub A$.
Since $A\in [X]^{<\kappa}$ was chosen arbitrarily,
by Lemma~\ref{lem:general_way_of_proof} there are $U_{n,m} \in \cU_{n,m}$ for all $n,m \in \w$ such that 
$\sset{U_{n,m}}{n,m \in \w} \in \Omega(\michk{X \cup \Fin}{k} \x Y)$.

(2)
The product space $(X \cup \Fin)^k \x Y$ is a continuous image of $\michk{X \cup \Fin}{k} \x Y$.
Apply~(1) and the fact that $\sgom$ is preserved by continuous mappings.
\epf

\section{Applications to finite powers of scale sets}\label{sec:applications}

In the theory of combinatorial covering properties, the principal objects of study are 
sets of reals.
One may ask why we considered the Michael topology, rather than the standard topology,
in the previous section.
The reason lies in the following simple but nonobvious corollary of Theorem~\ref{thm:kappa_fin_unb_x_sbgg_is_sgom}, given in Proposition~\ref{prp:mainmich}, which often allows us to consider products of spaces with combinatorial structure rather than only powers of such spaces.
Using the Michael topology on $\PN$ can also simplify proofs concerning sets of reals.

\bprp\label{prp:mainmich}
Let $\kappa$ be a regular cardinal with $\kappa \leq \fd$, let $k$ be a positive natural number, let $X_0,\dotsc, X_k\sub\roth$ be sets such that the set $X=\Un_{i<k}X_i$ is a $\kappa$-fin-unbounded set, and let $Y$ be a set satisfying $\sbgg$.
Then the following assertions hold.
\begin{enumerate}
\item 
The product space
\[
\prod_{i<k}\mich{X_i \cup \Fin} \x Y
\]
satisfies $\sgom$.
\item The product space 
\[
\prod_{i<k}(X_i \cup \Fin) \x Y
\]
satisfies $\sgom$.
\end{enumerate}
\eprp

\bpf
(1) 
By Theorem~\ref{thm:kappa_fin_unb_x_sbgg_is_sgom}, the product space 
\[
\prod_{i<k}\mich{X\cup \Fin}\x Y
\]
satisfies $\sgw$.
Then the product space $\prod_{i<k}\mich{X_i \cup \Fin} \x Y$ satisfies $\sgw$ because it is a closed subspace of the above product.

(2) The product space $\prod_{i<k}(X_i \cup \Fin) \x Y$ is a continuous image of $\prod_{i<k}\mich{X_i \cup \Fin} \x Y$.
Apply~(1) and the fact that $\sgom$ is preserved by continuous mappings.
\epf 

As we have already seen, the proof of Theorem~\ref{thm:kappa_fin_unb_x_sbgg_is_sgom} is highly nontrivial and technical.
Although the same method could be used to prove Proposition~\ref{prp:mainmich} directly, the direct proof would be more complicated because the factors may differ.
In the remainder of this section, we discuss applications of Proposition~\ref{prp:mainmich}.
To this end, we need examples of sets whose unions are $\kappa$-fin-unbounded sets for some regular cardinal $\kappa$.

\bprp
Let $\kappa$ be an infinite cardinal number and $X \sub \roth$ be a set with $\card{X} \geq \kappa$.
Then the following assertions are equivalent.
\begin{enumerate}
	\item The set $X$ is $\kappa$-fin-unbounded.
	\item For each $d \in \roth$
	there exists $S \sub X$ with $\card{S}<\kappa$ such that for each finite nonempty set $F \sub X \sm S$ we have $d\le^\infty \Un F$.
	\item For each $d \in \roth$
	there exists $S \sub X$ with $\card{S}<\kappa$ such that for each finite nonempty set $F \sub X \sm S$ we have $d\le^\infty \min F$.
\end{enumerate}
\eprp
\bpf
(1) $\Rightarrow$ (2)
Fix $d \in \roth$.
Let $\tilde{d} \in \roth$ be a function such that 
\[
	\tilde{d}(k) := \begin{cases}
		d(0),& \text{if } k=0, \\
		d(\tilde{d}(k-1)+1)+1, & \text{if } k>0.
	\end{cases}
\]
Then there exists a set $S \sub X$ with $\card{S} < \kappa$
such that for each finite nonempty set $F \sub X \sm S$ the set $\Un F$ omits infinitely many intervals of $\tilde{d}$.
Let $F \sub X \sm S$ be a finite nonempty set.
Then the set 
\[
I \coloneqq \sset{k \in \w}{[\tilde{d}(k), \tilde{d}(k+1)) \cap \Un F = \emptyset}
\]
is infinite and for every $k \in I$, we have 
\[
\Un F(\tilde{d}(k)+1) \geq \tilde{d}(k+1)>
d(\tilde{d}(k)+1).
\] 
Since $\sset{\tilde{d}(k)+1}{k \in I}$ is infinite, we have $d<^\infty \Un F$.

(2) $\Rightarrow$ (3) 
It follows from the fact that $\Un F \leq^\ast \min F$ for any finite nonempty set $F \sub X$.

(3) $\Rightarrow$ (1)
Let $d \in \roth$. There exists $t \in \roth$ such that for each natural $k \ge 1$ the interval $[t(k-1), t(k))$ contains $k^2+1$ consecutive intervals of $d$.
Then there exists a set $S \sub X$ with $\card{S}<\kappa$ such that for each nonempty finite set $F \sub X \sm S$ we have $t \le^\infty \min F$.
Thus, the set 
\[
J \coloneqq \sset{k \in \w}{t(k) \le x(k) \text{ for all } x \in F}
\]
is infinite.
Take $k \in J$ such that $k>\card{F}$. Then $\card{x \cap [t(k-1), t(k))} \leq k$ and thus $\card{\Un F \cap [t(k-1), t(k))}<k^2$.
Since $[t(k-1), t(k))$ contains more than $k^2$ consecutive intervals of $d$, the set $\Un F$ omits some interval of $d$ in $[t(k-1), t(k))$. Since $J$ is infinite, the set $\Un F$ omits infinitely many intervals of $d$. 
\epf

For $a\in \PN$, define $a\comp:=\w\sm a$.
Let $F$ be a semifilter. 
By $F^+$ we denote the set $\set{a \in \roth}{a\comp \notin F}$.
There is a notion in the literature that generalizes the $\bof(F)$-scale mentioned in the introduction.
A set $X\sub\roth$ is an \emph{$F$-scale}~\cite[Definition 4.1]{ST} if $\card{X}\geq \bof(F)$ and for any function $d\in\roth$ there exist a function $c\in\roth$ and a set $S\sub X$ of cardinality smaller than $\bof(F)$ such that 
\[
d \leq_{F^+} c \leq_F x
\]
for all functions $x\in X\sm S$.
It is easy to see that 
if $X=\{x_\alpha:\alpha<\mathfrak b(F)\}$ is 
a $\mathfrak b(F)$-scale, then it is an $F$-scale.
Indeed, given $d$, let $\gamma<\mathfrak b(F)$ be such that $x_\gamma\not\leq_F d$, which is equivalent to $d\leq_{F^+}x_\gamma$. Then $c=x_\gamma $ and $S=\{x_\xi:\xi\leq\gamma\}$ are as required, i.e., $c\leq_F x$ for all $x\in X\sm S$.

\bprp\label{lem:F-scale_is_b(f)_unbdd}
Let $F$ be a filter.
Each $F$-scale is a 
$\fb(F)$-fin-unbounded set. 
\eprp
\bpf
Take any $d \in \roth$. There exists $c \in \roth$ and a set 
$S \sub X$ with $\card{S}<\fb(F)$ such that 
\[
d \leq_{F^+} c \leq_F x
\]
for all functions $x \in X \sm S$. Now take any finite set $F_X \sub X \sm S$. Then $c \leq_F \min F_X$. Since $d \leq_{F^+} c$, it follows that $d \leq^\infty \min F_X$. Since $\card{X} \geq \fb(F)$, the set $X$ is a $\fb(F)$-fin-unbounded {set}.
\epf

\bcor
Let $F$ be a filter, let $X$ be an $F$-scale, let $k$ be a natural number, and let $Y$ be a set satisfying $\sbgg$.
Then the product space $\michk{X \cup \Fin}{k} \x Y$ is $\sgom$.
\ecor

In general, we cannot prove Proposition~\ref{lem:F-scale_is_b(f)_unbdd} for semifilters. To explain this we need 
to introduce the notion of \emph{bi-$\fd$-unbounded} set. For $a, b \in \roth$ we say that 
$a \le b$ if $a(n) \le b(n)$ for each $n \in \w$. 
Let $\kappa$ be an infinite cardinal number.
We say that $X \sub \roth$ is \emph{$\kappa$-unbounded} if
$\card{X} \geq \kappa$ and for each $d \in \roth$ the cardinality of the set $\sset{x \in X}{x \le d}$ is smaller than
$\kappa$. 
If $\cof(\kappa)>\w$ 
then $X \sub \roth$ is 
$\kappa$-unbounded if and only if $\card{X} \ge \kappa$ and $\card{\set{a \in X}{a \les c}} < \kappa$ for each $c \in \roth$.
We say that $X \sub \roth$ is \emph{bi-$\fd$-unbounded} if the sets $X$ and $\sset{x\comp}{x \in X} \sub \roth$
are $\fd$-unbounded. Szewczak and Tsaban proved~\cite[Theorem 3.3]{ST} that bi-$\fd$-unbounded subsets of 
$\roth$ exist if and only if $\fd \le \fr$.
\bexm
Let $F:= \roth$ be the full semifilter and $X' \sub \roth$ be a bi-$\fd$-unbounded set.
Then $X \coloneqq X' \cup \sset{x\comp}{x \in X'}$ is an $F$-scale that is not a $\fb(F)$-fin-unbounded set.
\eexm

\bpf
First, we show that $X$ is an $F$-scale. Note that $\fb(F) = \fd$ so $\card{X} \geq \fb(F)$. 
Take any $d \in \roth$.
We may assume that $d\neq \w$.
For every $k \in \w$ define $d_k \in \roth$ by $d_k(n) \coloneqq d(k+n)$ for all
natural numbers $n$.
Let $D \coloneqq \sset{d_k \in \roth}{k \in \w}$ and let
\[
S \coloneqq \sset{x \in X'}{x \le d \text{ for some } d \in D} \cup \sset{x\comp}{x \in X', x\comp \le d \text{ for some } d \in D}.
\]
Then $\card{S}<\fd$ and for all $x \in X \sm S$ we have $x \not\le d_k$ for all $k \in \w$.
It follows that $d \le^\infty x$ (equivalently $d \le_F x$).
Indeed, take any natural number $n$. Then there exists a natural $k>n$ such that $x \not \le d_k$. Hence, there exists a natural number
$m$, such that 
\[
d(k+m)=d_k(m) < x(m) \le x(k+m).
\]
Since $k+m>n$, it follows that $d \le^\infty x$. 
Thus, $X$ is an $F$-scale.

Now suppose that $X$ is a $\fb(F)$-fin-unbounded set and let $d\in\roth\sm\{\w\}$.
Then there exists $S \sub X$ with $\card{S}<\fd$ such that for all nonempty finite sets $H \sub X \sm S$ we have $d \le^\infty \Un H$.
Since $\card{S}<\fd$, it follows that there exists $x \in X \sm S$ such that $x\comp \in X \sm S$. 
Then $\{x, x\comp\} \sub X \sm S$ and thus $d \le^\infty x \cup x\comp = \w$, a contradiction. 
\epf

\bprp
Let $U$ be an ultrafilter. A union of fewer than $\fb(U)$ $U$-scales is also a $U$-scale.
In particular, it is a $\fb(U)$-fin-unbounded set.
\eprp
\bpf
Let $X=\Un_{\alpha < \kappa} X_\alpha$ where the sets $X_\alpha$ are $U$-scales and $\kappa < \fb(U)$. Take any function $d \in \roth$. Then for each 
$\alpha<\kappa$ there are a function $c_\alpha \in \roth$ and a set $S_\alpha \sub X_\alpha$ of cardinality smaller than $\bof(U)$ such that 
\[
d \leq_{U^+} c_\alpha \leq_U x
\]
for all functions $x \in X_\alpha \sm S_\alpha$. Let 
\[
S \coloneqq \Un_{\alpha<\kappa} S_\alpha.
\]
Since $U^+=U$, it follows that $d \leq_U c_\alpha \leq_U x$ for all 
$x \in X_\alpha \sm S$ and $\alpha<\fb(U)$.
Then $d \leq_U x$ for all $x \in X \sm S$.
Since $\bof(U)$ is regular, we have $\card{S}<\bof(U)$.
\epf
\bcor
Let $U$ be an ultrafilter, let $k$ be a natural number, let $X_0,\dotsc,X_k$ be $U$-scales, and let $Y$ be a set satisfying $\sbgg$.
Then the product space $\prod_{i\leq k}\mich{X_i \cup \Fin} \x Y$ satisfies $\sgom$.
\ecor

\bcor\label{cor:product_of_U-scales}
Let $U$ be an ultrafilter, let $k$ be a natural number, let  $X_0,\dotsc,X_k$ be $U$-scales, and let $Y$ be a set satisfying $\sbgg$.
Then the product space $\prod_{i\leq k} (X_i \cup \Fin) \x Y$ satisfies $\sgom$.
\ecor

\bthm[Talagrand~\cite{T1980}]\label{thm:Talagrand_char}
Let $F$ be a semifilter.
The following assertions are equivalent.
\be 
\item The semifilter $F$ is meager;
\item There is a function $h\in\roth$ such that for every $x \in F$ and for all but finitely many natural numbers $n$
 we have $x \cap [h(n), h(n+1)) \neq \emptyset$.
\ee
\ethm

\bprp
Let $\{F_\xi :\xi<\kappa\}$   be a family of  meager semifilters for some $\kappa<\fb$, and $X_\xi$ be a $F_\xi$-scale for every $\xi\in \kappa$. Then $X=\bigcup_{\xi\in\kappa}X_\xi$ is 
$\fb$-fin-unbounded. 
\eprp
\bpf
Suppose, contrary to our claim, that there exists $d\in\roth$ such that for every
$S\in [X]^{<\fb}$, there exists a finite set $K\in [X\setminus S]^{<\w}$ such that $\Un K\in L$, where 
\[
L \coloneqq \sset{c \in \roth}{c \text{ meets all but finitely many intervals of }d}.
\]
Note that $L$ is a $\sigma$-compact subset of $\roth$.
Furthermore, since $F_\xi$ is meager for all $\xi<\kappa$, by~\cite[Corollary~2.27]{sfh} we have $\fb(F_\xi) = \fb$ for all $\xi<\kappa$. 
Therefore 
$X_\xi$ is $\fb$-unbounded for all $\xi<\kappa$, and thus $X$ is $\fb$-unbounded.
Define 
${\downarrow}c \coloneqq \sset{c' \in \roth}{c'\les c}$ for all $c\in\roth$.
For a positive $n$, let 
\[
C_n \coloneqq \set{c \in \roth}{\text{there exists an $n$-element set } K_c \sub X \sm {\downarrow} c, \text{ such that } \Un K_c \in L}.
\]
Observe that $\Un_{n =1}^\infty C_n = \roth$.
Indeed, since $X$ is $\fb$-unbounded, it follows that for every $c\in\roth$, the cardinality of the set 
$X\cap{\downarrow} c$ is smaller than $\fb$.
Thus, there exists a finite set $K_c \sub X\setminus{\downarrow} c$,
such that $\Un K_c\in L$.
Thus, $c \in \Un_{n =1}^\infty C_n$.

Since the union of countably many non-dominating sets is non-dominating, there exists $n>0$ such that $C_n$
is dominating. Let $n_0$ be the smallest natural number such that $C_{n_0}$ is dominating. 
In particular,
there exists $c_0\in\roth$ such that 
$K\cap {\downarrow} c_0\neq\emptyset$ for every $K\in [X]^{<n_0}$
with $\Un K\in L$.

Define $\varphi \colon ((X\cup \Fin)\setminus{\downarrow} c_0)^{n_0}\to \PN$ by $\varphi(y)=\bigcup_{i\in n_0}y(i)$ and note that $\varphi$ is continuous.
The set $(X\cup \Fin)\setminus{\downarrow} c_0$ is productively Hurewicz in the realm of metrizable, separable spaces. Indeed, $(X_\xi \cup \Fin) \sm {\downarrow} c_0$ is an $F_\xi$-scale, thus by~\cite[Theorem~5.4]{ST} 
the space $(X_\xi \cup \Fin) \sm {\downarrow} c_0$ is productively Hurewicz in the realm of metrizable, separable spaces for all $\xi<\kappa$.
Note that for each $f \in \kappa^{n_0}$ the space
\[
Z_f \coloneqq \prod_{i<n_0}((X_{f(i)} \cup \Fin) \sm {\downarrow} c_0)
\]
is Hurewicz. Observe that 
\[
\big((X \cup \Fin)\sm {\downarrow} c_0\big)^{n_0} = \Un_{f \in \kappa^{n_0}} Z_f.
\]
Since the Hurewicz property is preserved by unions of families of size smaller than $\fb$, it follows that $((X \cup \Fin) \sm {\downarrow} c_0)^{n_0}$ is Hurewicz.
Consequently, the subspace $T:=\varphi^{-1}[L]$ is Hurewicz, being an $F_\sigma$-subset of $((X\cup \Fin)\setminus{\downarrow} c_0)^{n_0}$. 
We claim that $T\sub (X\setminus{\downarrow} c_0)^{n_0}$. Indeed,
if $y\in ((X\cup \Fin)\setminus{\downarrow} c_0)^{n_0}$ is such that $y(i)\in\Fin$
for some $i\in n_0$, then $\varphi(y)=^*\Un K$
for some $K\in [X\setminus{\downarrow} c_0]^{<n_0}$.
Hence 
$\Un K=^*\varphi(y)\not\in L$ because $K\cap{\downarrow} c_0=\emptyset$.
Since $T$ is an $F_\sigma$ subset of $((X\cup\Fin)\setminus{\downarrow} c_0)^{n_0}\sub(\roth)^{n_0}$ and $((X\cup\Fin)\setminus{\downarrow} c_0)^{n_0}\sub(\roth)^{n_0}$ is Hurewicz,  $T$ is also Hurewicz. It follows that $T$ is a Hurewicz subset of 
$(X\setminus{\downarrow} c_0)^{n_0}\sub(\roth)^{n_0}$, and hence
there exists $c\in\roth $ such that
$T\sub({\downarrow} c)^{n_0}$. Since $C_{n_0}$ is dominating, without loss of generality we may assume that
 $c_0\leq^*c$ and $c\in C_{n_0}$.  By the definition of $C_{n_0}$
there exists  $K_c \sub X \sm {\downarrow}c$ of size  $|K_c|=n_0$
with $\Un K_c \in L$. Let $t \colon n_0\to K_c$ be a bijection. 
Then $t \in T$ and $t(i) \notin {\downarrow} c$ for all $i \in n_0$, a contradiction.
\epf

\bcor
Let $\sset{F_i}{i \leq  k}$ be a finite family of meager semifilters and let $X_i$ be an $F_i$-scale for every $i \leq k$.
Let also $Y$ be a set satisfying $\sbgg$.
Then the product space $\prod_{i \leq k}\mich{X_i \cup \Fin} \x Y$ satisfies $\sgom$.
In particular, it holds if $X_0,\dotsc,X_{k}$ are $\fb$-scales.
\ecor

\section{$\gamma$-sets}

In this section, we consider properties from the third row of the Scheepers diagram.
Let $\Omctbl$ be the class of all countable $\w$-covers of spaces.
Our main goal is to show the following result.
\bthm\label{thm:gamma}
Let $\kappa$ be an infinite regular cardinal with $\kappa \le \cov(\cM)$, let $X$ be a $\kappa$-fin-unbounded set, let $k$ be a natural number, and let $Y \sub \PN$ be a set with the $\gamma$-property.
\be
\item The space $\michk{X \cup \Fin}{k} \x Y$ satisfies $\sone(\Omctbl, \Om)$.
\item The space $(X \cup \Fin)^k \x Y$ satisfies $\sone(\Om, \Om)$.
\ee
\ethm

A straightforward application of Theorem~\ref{thm:gamma}, given in Proposition~\ref{prp:gamma}, generalizes a result of Szewczak and W{\l}udecka~\cite[Theorem 5.6]{unbddtower}.
To obtain this result, we require the following theorems and auxiliary results, some of which are derived from the results of 
Szewczak--W{\l}udecka~\cite{unbddtower}.
\bthm\cite[Theorem 2]{GN1982}\label{thm:Gerlits-Nagy}
Let $X$ be a a subspace of a space $Y$. Then the following implications hold
\[
\text{(1)}\Rightarrow \text{(2)}\Rightarrow \text{(3)},
\]
where
\be 
\itm $X$ has the $\gamma$-property;
\itm $X$ has the property $\sone(\Omprim, \Gammaprim)$; 
\itm For each $\cU \in \Omegaprim(X)$ there is a family $\cV \sub \cU$ such that $\cV \in \Gammaprim(X)$.
\ee
\ethm

\blem\label{lem:Quasi_sone_omega_omega}
Let $A$ be a set with $\card{A}<\cov(\cM)$, let $k$ be a natural number, let $Y$ be a set with the $\gamma$-property, and let $\cU_0, \cU_1, \dotsc \in \Omctblprim( {\michk{\Fin \cup A}{k}}\x Y)$.
Then there exists a family $\set{U_n}{n \in \w}$ with $U_n \in \cU_n$ such that 
for each finite set $H \sub A$ there is a family $\cU_H \sub \set{U_n}{n \in \w}$ such that $\cU_H \in \Gammaprim({\michk{\Fin \cup H }{k}}\x Y)$.
In particular, we have 
\[
\sset{U_n}{n\in\w}\in \Omprim({\michk{\Fin \cup H }{k}}\x Y).
\]
\elem
\bpf
Take any finite set $H \sub A$ and let $\cU_n = \set{U_{n,m}}{m \in \w} \in \Omctblprim({\michk{\Fin \cup A}{k}}\x Y)$.
Since ${\michk{\Fin \cup H}{k}} \x Y$ has the $\gamma$-property, {by Theorem~\ref{thm:Gerlits-Nagy},} there are sets $U_{n, f_H(n)} \in \cU_n$ for each $n$ such that $\set{U_{n, f_H(n)}}{n \in \w} \in \Gammaprim( {\michk{\Fin \cup H}{k}} \x Y)$.
Note that $\card{\set{f_H \in \w^\w}{H \in [A]^{<\w}}}<\cov(\cM)$; thus there exists $f \in \w^\w$ such that 
for each finite set $H \sub A$ the set \[N_H \coloneqq \set{n \in \w}{f_H(n)=f(n)}\] is infinite. 
Note that
\[\cU_H \coloneqq \set{U_{n, f_H(n)}}{n \in N_H} \sub \set{U_{n ,f(n)}}{n \in \w}\] and $\cU_H \in {\Gammaprim(\michk{H \cup \Fin}{k} \x Y)}$.
\epf
\blem\label{lem:Galvin-Miller}\cite[Lemma 1.2]{GM1984}
Let $\cU$ be a family of open sets in $\PNM$ such that for each finite set $F\sub\Fin$, there exists infinitely many $U\in\cU$ such that $F\sub U$.
Then there exist a function $a \in \roth$ and pairwise distinct sets $U_0, U_1, \dotsc \in \cU$ such that for each set $x \in \roth$ and all natural 
numbers $n$:
\[
 \text{If } x \cap [a(n), a(n+1)) = \emptyset, \text{ then } x \in\Intc(U_n).
\]
\elem

The original formulation of Lemma~\ref{lem:Galvin-Miller} is stated for the space $\PN$ instead of $\PNM$.
However, the same proof works for our modification as well.

Let $\cU$ and $\cV$ be covers of the space $X$. We say that $\cV$ \emph{refines} $\cU$ if for each 
$V \in \cV$, there exists $U \in \cU$ such that $V \sub U$.
Let $X$ and $Y$ be topological spaces. Their disjoint sum $X\sqcup Y$ is the set 
\[
(X\x\{0\})\cup(Y\x\{1\})
\]
equipped with the topology in which a set $U\sub X\sqcup Y$ is open if and only if $U\cap(X\x\{0\})$ is open in $X\x\{0\}$ and $U\cap(Y\x\{1\})$ is open in $Y\x\{1\}$.
Although the spaces $\PNM$ and $\PN$ have the same underlying set, we regard them as disjoint copies of this set when considering $\PNM\sqcup\PN$. To simplify notation, we suppress the second coordinate throughout.
Thus, we write elements of $\PNM\sqcup\PN$ simply as elements of $\PNM$ or $\PN$, with the ambient copy determined by the space to which the element is declared to belong.
\blem\label{lem:refinement2}\cite[Lemma 5.2]{unbddtower} 
Let $Y \sub \Pof(\w)$, $H \sub \roth$ be a finite set and $k$ be a natural number. 
Let $\cU \in \Omctblprim(\michk{\Fin \cup H}{k} \x Y)$ be a family of open sets in $\michk{\PN}{k}\x \PN$.
There is a family $\cV \in \Omctblprim(\mich{\Fin \cup H} \sqcup Y)$ of open sets in $\PNM \sqcup \PN$ such that the family
\[
\set{V^{k+1} \cap (\michk{\PN}{k} \x Y)}{V \in \cV} \in \Omctblprim(\michk{\Fin \cup H}{k} \x Y)
\]
refines the family $\cU$.
\elem
\bpf
Let $\cU \in \Omctblprim(\michk{\Fin \cup H}{k} \x Y)$ be a family of open sets in $\michk{\PN}{k}\x \PN$.
Let $\cB$ be a countable base of $\PN$ that is closed with respect to finite unions and let $\cB_\mathrm{M}:=\sset{B\cup H}{B\in\cB}$.
Let $F_0 \sub \Fin\cup H$ and let $F_1 \sub Y$ be non-empty finite sets such that $F_0\cap\Fin\neq\emptyset$. 
Since $\cU \in \Omctblprim(\michk{\Fin \cup H}{k} \x Y)$, there is a set $U\in\cU$ such that $F_0^k \x F_1 \sub U$.
There are sets $B_0\in \cB_\mathrm{M}$ and $B_1\in \cB$ such that
\[
F_0^k\x F_1\sub B_0^k\x B_1\sub U.
\]
Note that there are infinitely many pairs $B_0,B_1$ with the above property.
It follows that the family
\[
\smallmedset{B_0^k\x B_1}{B_0\in\cB_M, B_1\in \cB, B_0^k\x B_1\sub U\text{ for some  }U\in\cU}
\]
is in $\Omctblprim(\michk{\Fin \cup H}{k} \x Y)$ and refines $\cU$.
For every $B_0\in \cB_M$ and $B_1\in\cB$, we have
\[
B_0\x B_1=(B_0\sqcup B_1)^{k+1}\cap (\PN_{\mathrm{M}}^k \x \PN).
\]
It follows that the family
\[
\smallmedset{B_0^k\sqcup B_1}{B_0\in\cB_M, B_1\in \cB, B_0^k\x B_1\sub U\text{ for some  }U\in\cU}
\]
has desired properties.
\epf

\blem\label{lem:Kappa-fin-unbdd_And_Gamma_set} 
Let $\kappa$ be an infinite cardinal number and let $H \sub \roth$ be a finite set. Let $Y \sub \PN$ be a set satisfying the $\gamma$-property and let $X \sub \roth$ be a $\kappa$-fin-unbounded set.
For each family $\cU \in \Omctblprim(\mich{\Fin \cup H} \sqcup Y)$ of open sets in $\PN_\mathrm{M} \sqcup \PN$ there is a set $S \sub X$ with 
$\card{S}<\kappa$ such that 
\[
\cU \in \Omctblprim(\mich{\Fin\cup(X \sm S) \cup H} \sqcup Y).
\]
\elem

\bpf 
Let $\cU \in \Omctblprim(\mich{\Fin \cup H} \sqcup Y)$ be a family of open sets in $\PN_\mathrm{M} \sqcup \PN$. 
Note that $\mich{\Fin \cup H} \sqcup Y$ has the $\gamma$-property.
By Theorem~\ref{thm:Gerlits-Nagy}, there exists $\cU' \sub \cU$ such that 
$\cU' \in \Gammaprim(\mich{\Fin \cup H} \sqcup Y)$. 
Each set $U\in\cU'$ is of the form $U=U^{(0)}\sqcup U^{(1)}$, where $U^{(0)}\sub\PN_\mathrm{M}$ and $U^{(1)}\sub\PN$ are open.

Apply Lemma~\ref{lem:Galvin-Miller} to the family $\sset{U^{(0)}}{U\in\cU'}$. Then there exist a function $a \in \roth$ and sets 
$U_0, U_1, \dotsc \in \cU$ such that, for every $x \in \PNM\cap \roth$ and every $n \in \w$,
\[
 \text{if } x \cap [a(n),a(n+1)) = \emptyset, \text{ then } x \in \Intc(U_n^{(0)}).
\]
Since $X$ is a $\kappa$-fin-unbounded set, there exists a set $S \sub X$ with $\card{S}<\kappa$ such that for all finite
sets $F \sub X \sm S$ we have $\Un F$ omits infinitely many intervals of $a$. Thus, for all finite 
sets $F \sub X \sm S$ we have $F \sub \Intc(U_n^{(0)})$ for infinitely many $n \in \w$.

Since $\cU' \in \Gamma(\mich{\Fin \cup H} \sqcup Y)$, for every finite set 
\[
K\sub(\Fin\cup H)\sqcup Y
\]
we have $K\sub U_n$ for all but finitely many $n$.
It follows that for every finite set $F\sub X\sm S$ and every finite set 
$K\sub(\Fin\cup H)\sqcup Y$, there {exist infinitely many natural numbers} $n$ such that $F\cup K\sub U_n$.
Thus,
\[
\cU' \in \Omctblprim(\mich{\Fin \cup H \cup X \sm S} \sqcup Y).
\]
But then
\[
\cU \in \Omctblprim(\mich{\Fin \cup H \cup X \sm S} \sqcup Y).\qedhere
\]
\epf

\blem\label{lem:Kappa-fin-unbdd_And_Gamma_Set_Powers}
Let $\kappa$ be a cardinal number, let $X \sub \roth$ be a $\kappa$-fin-unbounded set, let $Y \sub \PN$ be a set satisfying the $\gamma$-property, let $H \sub \roth$ be a finite set, and let $k$ be a natural number.
For each family 
$
\cU \in \Omctblprim(\michk{\Fin \cup H}{k} \x Y)
$
of open sets in $\michk{\PN}{k}\x\PN$, there exists a set $S \sub X$ with 
$\card{S}<\kappa$ such that {each finite set $F \sub \michk{\Fin \cup H \cup X \sm S}{k} \x Y$ is contained in some member of $\cU$.}
\elem
\bpf
By Lemma~\ref{lem:refinement2}, there is a family $\cV \in \Omctblprim(\mich{\Fin \cup H} \sqcup Y)$
such that the family
\[
\set{V^{k+1} \cap (\michk{\PN}{k} \x {\PN})}{V \in \cV}
\]
refines the family $\cU$.
By Lemma~\ref{lem:Kappa-fin-unbdd_And_Gamma_set}, there exists a set $S \sub X$ with
$\card{S}<\kappa$ such that
$
\cV \in \Omctblprim(\mich{\Fin \cup H \cup X \sm S} \sqcup Y).
$
Take any finite sets
$F_X \sub \Fin \cup H \cup X \sm S$ and $F_Y \sub Y$.
{Then} there {exist infinitely many sets} $V\in\cV$ such that $F_X\cup F_Y\sub V$.
Then
\[
F_X^k \x F_Y
\sub
V^{k+1} \cap (\michk{\PN}{k} \x {\PN})\sub U
\]
for some set $U\in\cU$.
It follows that {each finite set $F \sub \michk{\Fin \cup H \cup X \sm S}{k} \x Y$ is contained in some member of $\cU$.}
\epf

\bpf[{Proof of Theorem~\ref{thm:gamma}}]
 (1)
 Let $\cU_{n,m} \in \Omctbl(\michk{X \cup \Fin}{k} \x Y)$ for all natural numbers $n,m$.
We prove that the families $\cU_{n,m}$ satisfy the assumptions of Lemmma~\ref{lem:general_way_of_proof}:
 Take any $n \in \w$ and any set $A \sub X$ with $\card{A}<\kappa$.
Since $\card{A}<\cov(\cM)$, Lemma~\ref{lem:Quasi_sone_omega_omega} yields $U_{n,m}(A) \in \cU_{n,m}$ such that 
for every finite set $H \sub A$, we have 
\[
\set{U_{n,m}(A)}{m \in \w}\in \Omprim(\michk{\Fin\cup H}{k} \x Y).
\]
By Lemma~\ref{lem:Kappa-fin-unbdd_And_Gamma_Set_Powers}, there are sets $S_H \sub X$ with 
$\card{S_H}<\kappa$ such that 
\[
\set{U_{n,m}(A)}{m \in \w} \in \Omegaprim(\michk{\Fin \cup H \cup X \sm S_H}{k} \x Y)
\]
for all finite sets $H\sub A$.

By Lemma~\ref{lem:general_way_of_proof}, there are sets $U_{n,m} \in \mathcal U_{n,m}$ such that 
\[
\set{U_{n,m}}{n,m \in \w} \in \Omega(\michk{X \cup \Fin}{k} \x Y).
\]

(2)
The product space $(X \cup \Fin)^k \x Y$ is a continuous image of $\michk{X \cup \Fin}{k} \x Y$.
Apply~(1), the fact that $\sone(\Omctbl, \Om)$ is preserved by continuous mappings, and the properties $\sone(\Omctbl, \Om)$ and $\sww$ are equivalent in the class of sets of reals.
\epf

\bprp\label{prp:gamma}
Let $\kappa$ be an infinite regular cardinal with $\kappa \leq \cov(\cM)$, let $k$ be a positive natural number, let $X_0,\dotsc, X_{k-1}\sub\roth$ be sets such that the set $X=\Un_{i<k}X_i$ is a $\kappa$-fin-unbounded set, and let $Y \sub \PN$ be a set satisfying the $\gamma$-property.
\be 
\item The product space
\[
\prod_{i<k}\mich{X_i \cup \Fin} \x Y
\]
satisfies $\sone(\Omctbl, \Omega)$.
\item The product space
\[
\prod_{i<k} (X_i \cup \Fin) \x Y
\]
satisfies $\sone(\Omega, \Omega)$.
\ee
\eprp
\bpf
(1) 
By Theorem~\ref{thm:gamma}, the product space 
\[
\prod_{i<k}\mich{X\cup \Fin}\x Y
\]
satisfies $\sone(\Omctbl, \Omega)$.
Then the product space $\prod_{i<k}\mich{X_i \cup \Fin} \x Y$ satisfies $\sone(\Omctbl, \Omega)$ because it is a closed subspace of the above product.

(2) 
The product space $\prod_{i<k}(X_i \cup \Fin) \x Y$ is a continuous image of $\prod_{i<k}\mich{X_i \cup \Fin} \x Y$.
The property $\sone(\Omctbl, \Om)$ is preserved by continuous images. 
Apply~(1) and the fact that the properties $\sone(\Omctbl, \Omega)$ and $\sone(\Omega, \Omega)$ are equivalent in the realm of sets of reals.
\epf 
\bcor
Let $\kappa$ be an infinite regular cardinal with $\kappa \le \cov(\cM)$ and let
$X$ be a $\kappa$-fin-unbounded set. Let $k$ be a natural number and let $Y \sub \PN$ be a set satisfying the $\gamma$-property.
Then $(X \cup \Fin)^k \x Y$ is $\sone(\Omega, \Omega)$.
\ecor
This corollary generalizes the result of Szewczak, Tsaban, and Zdomskyy~\cite[Lemma 2.21]{stz} that was mentioned in the introduction to this paper, since the one-point space has the $\gamma$-property. 
It is noteworthy that their approach uses the definition of $U$-Menger space and its topological properties, while our approach is mainly combinatorial and is therefore much more flexible. It is also a generalization of results by Miller, Tsaban, and Zdomskyy~\cite[Corollary 6.9~(2)]{SPMProd} combined with a result of Ko\v{c}inac and Scheepers~\cite[Theorem 19]{CoC7}. Their results also use advanced topological tools rather than combinatorial ones.
\bcor
Let $\kappa$ be an infinite regular cardinal with $\kappa \le \cov(\cM)$ and let 
$X$ be a $\kappa$-fin-unbounded set. Then $X \cup \Fin$ is $\sone(\Omega, \Omega)$.
\ecor

\bcor 
Let $\kappa$ be an infinite regular cardinal with $\kappa \le \cov(\cM)$, let
$X$ be a $\kappa$-fin-unbounded set, and let $k$ be a natural number.
Then $(X \cup \Fin)^k$ is Rothberger.
\ecor

\section{Counterexamples}\label{sec:counterexamples}
In this section, we present two counterexamples.
Our aim is to show that in Theorem~\ref{thm:kappa_fin_unb_x_sbgg_is_sgom} the property $\sbgg$ cannot be replaced by $\sbomom$, even if $X$ is a $\fb$-scale. The first example is constructed using ideas from papers by Szewczak and Wi\'{s}niewski~\cite{luzincomb} and Bartoszy\'{n}ski, Shelah, and Tsaban~\cite{BST2003}.
The second example is constructed using Cohen forcing.

\bthm\label{thm:counter1}
Assume that $\cov(\cM)=\fc$ and $\fc$ is regular.
Then there exist ultrafilters $U$ and $\tilde{U}$, a $U$-scale $X$, and a $\tilde{U}$-scale $Y$, such that 
$X \cup \Fin$ and $Y\cup\Fin$ are $\sbomom$, but the product space $(X \cup \Fin) \times (Y\cup\Fin)$ is not Menger. In fact, even $X$ and $Y$ satisfy the stated property, and $X \x Y$ is not Menger. 
\ethm

By $[\w]^{\w, \w}$ we denote the family of all infinite and co-infinite subsets of $\w$.
\bprp\label{prop:comb_luzin}\cite[Lemma 2.4]{luzincomb}
Let $\cM'$ be a family of meager subsets of $\roth$ with $\card{\cM'} < \cov(\cM)$. 
Then for each family $Z \sub \roth$ with $\card{Z} < \cov(\cM)$ and each $d \in \roth$ 
there exist $x, y \in [\w]^{\w, \w} \sm \Un \cM'$ such that 
\[
z <^\infty x, z <^\infty y \text{ for each } z \in Z \text{ and } d \le^\ast \max\{x\comp, y\comp\}.
\]
\eprp
\blem\label{lem:comb_luzin_fin_mod}
Let $\cM', Z$, and $d$ be as in Proposition~\ref{prop:comb_luzin}. 
Then there are $x, y \in [\w]^{\w, \w}$ such that 
\[
z <^\infty x', z <^\infty y' \text{ for each } z \in Z,\quad d \le^\ast \max\{(x')\comp,(y')\comp\},\quad\text{ and }x',y'\notin \Un\cM'
\]
for all finite modification $x',y'$ of $x$ and $y$, respectively.
\elem
\bpf
We may assume that the set $\Un\cM'$ is closed under finite modifications of its elements.
Let $T=\{h(z):z\in Z\}$, where $h(z)(n)=z(2n)$ for all $n\in\w$.
Since $\card{T}<\cov(\cM)$ by Proposition~\ref{prop:comb_luzin}, 
there exist $x, y \in [\w]^{\w, \w} \sm \Un \cM'$ such that 
\[
t <^\infty x, t<^\infty y \text{ for each } t \in T \text{ and } h(d) \le^\ast \max\{x\comp, y\comp\}.
\]
We claim that $x,y$ are as required.
Indeed, suppose that $x',y'$ are finite modifications of $x,y$, respectively.
Then there are $k_0,l_0,k_1,l_1\in\bbZ$ such that
$x'(i)=x(i+k_0)$, $y'(i)=y(i+l_0)$, $(x')\comp(i)=x\comp(i+k_1)$
and $(y')\comp(i)=y\comp(i+l_1)$ for all but finitely many $i \in \w$. In particular,
\[
\max\{ (x')\comp(i),  (y')\comp(i) \}\geq  \max\{ x\comp(i+p),  y\comp(i+p) \}
\] 
where $p=\min\{k_1,l_1\}$, for all but finitely many $i \in \w$.

Now, if  $x(i)>h(z)(i)=z(2i)$ and $i>|k_0|$, then
\[  x'(2i)> x'(i-k_0)=x(i+k_0-k_0)=x(i)>z(2i), \]
and hence $x'(2i)>z(2i)$ for infinitely many $i$.
The same argument applies to $y$.

Finally, if $\max\{ x\comp(i),  y\comp(i) \}>h(d)(i)=d(2i)$ and $i>|p|$,
then
\begin{gather*}
 \max\{ (x')\comp(i+|p|),  (y')\comp(i+|p|) \}\geq 
\max\{ x\comp(i+p +|p| ),  y\comp(i+p +|p| ) \} \geq\\
\max\{ x\comp(i ),  y\comp(i ) \} \geq d(2i) \geq
d(i+|p|) , 
\end{gather*}
hence $ \max\{ (x')\comp,  (y')\comp \} \geq^* d $.

Since $\Un\cM'$ is closed under finite modifications of its elements and $x,y\notin \Un\cM'$, we have $x',y'\notin\Un\cM'$.
\epf

Let $\cU \in \Om_{\mathrm{Borel}}(X)$. We say that $\cU$ is $\w$-fat if for every natural $k$, every finite $F \sub X$, and every finite family of nonempty open sets $O_1, \dotsc, O_k$, there exists $U \in \cU$ such that $F \sub U$ and none of the sets $U \cap O_1, \dotsc, U \cap O_k$ is meager. Let $\Om_{\mathrm{Borel}}^{\mathrm{fat}}$ be the collection of all countable $\w$-fat Borel covers of $X$.
We use the following elementary properties; for proofs, see Tsaban~\cite{T2006}.
\blem\cite[Lemma~3.9]{T2006} 
Assume that $\cU$ is a countable collection of Borel sets. Then $\Un \cU$ is comeager if and only if for each nonempty basic open set $O$ there exists $U \in \cU$ such that $U \cap O$ is not meager.
\elem
\bcor\cite[Corollary~3.10]{T2006}\label{cor:omega_fat_corrolary1}
Assume that $\cU \in \Om_{\mathrm{Borel}}^{\mathrm{fat}}(X)$. Then: 
\be 
\itm For each finite set $F \sub X$ and each finite family $\cF$ of nonempty basic open sets, the set 
\[
 \Un \set{U \in \cU}{F \sub U \text { and for each } O \in \cF, U \cap O \notin \cM}
\]
is comeager.
\itm For each element $x$ in the intersection of all sets of this form, $\cU$ is an $\w$-fat cover of $X \cup \{x\}$.
\ee
\ecor 
\blem\cite[Lemma 10]{BST2003}\label{lem:sbfombfom}
If $Y \sub \roth$ has cardinality less than $\cov(\cM)$, then $Y$ satisfies $\sbomfomf$.
\elem
\bpf
Suppose that $\eseq{\cU}$ is a sequence of $\omega$-fat covers of $Y$, where $\cU_n = \set{U_{n,k}}{k \in \w}$.
Let $\cB$ be a countable base of $\roth$.
For each finite set $F \sub Y$ and each finite $B \sub \cB$, let $f_{F,B}(n) = \min \set{k}{F \sub U_{n,k} \text{ and } V \cap U_{n,k} \text{ is not meager for each } V \in B}$.
The set \[H_{F,B} = \set{x \in \w^\w}{\forall^\infty n\ x(n) \neq f_{F,B}(n)}\] is meager in $\w^\w$.
It suffices to show that sets $G_k = \set{x \in \w^\w}{\forall n > k\ x(n) \neq f_{F,B}(n)}$ are nowhere dense for each $k \in \w$.
Take any $k \in \w$ and $x \in G_k$. Take any basic neighborhood $O_x$ of $x$. We may assume that 
$O_x = \set{y \in \w^\w}{y(m) = x(m) \text{ for all } m=0, 1, \dotsc, n_0}$. Then there exists $k_0>\max\{k,n_0\}$ such that 
$x(k_0) \neq f_{F,B}(k_0)$. Then the set
\[
 \set{y \in \w^\w}{y(m) = x(m) \text{ for } m=0, 1, \dotsc, n_0 \text{ and } y(k_0) = f_{F,B}(k_0)}
\]
is contained in $O_x$ but its intersection with $G_k$ is empty; thus $G_k$ is nowhere dense. Since 
$H_{F, B} = \Un_{k \in \w} G_k$, it is a meager subset of $\w^\w$.

Let \[H= \Un_{F \in [Y]^{<\w}, B \in [\cB]^{<\w}} H_{F,B}.\] Since $\card{Y}<\cov(\cM)$, it follows that the set 
$\w^\w \sm H$ is nonempty. Observe that any $z \in \w^\w \sm H$ gives the desired selector.
\epf
\blem\cite[Lemma~3.12]{T2006}\label{lem:all_borel_is_fat}
Assume that $X$ is a set of reals such that for each nonempty basic open set $O$, the set $X \cap O$ is not meager. 
Then every countable Borel $\w$-cover $\cU$ of $X$ is an $\w$-fat cover of $X$.
\elem
Let $Y\sub \w^\w$. We define $\maxfin(Y)$ to be the set of all $a \in \w^\w$ such that 
$a(k)=\max_{x \in F}x(k)$ for some finite set $F \sub Y$ and all $k \in \w$.
For $x,y\in [\omega]^\omega$ and a relation $R$ on $\omega$
we use $[x\, R\, y]$ to denote the set
$\{n\in\omega: x(n)R\,y(n)\}$.
\bpf[{Proof of Theorem~\ref{thm:counter1}}]
Let $\set{O_n}{n \in \w}$ be an enumeration of all nonempty basic open sets in $\roth$. Let $\set{\cF_m}{m \in \w}$ be an enumeration of all finite families
of nonempty basic open sets in $\roth$. 
Let $\sset{G_\alpha}{\alpha<\fc}$ be an enumeration of all dense $G_\delta$-sets in $\roth$.
Let $\set{d_\alpha}{\alpha<\fc}$ be a dominating set in $\roth$ and 
let $\seq{\seq{\cU_{\alpha, n}}{n \in \w}}{\alpha<\fc}$ enumerate all sequences of countable families of Borel sets.
We will define sets $F_\alpha, \tilde{F}_\alpha \sub \roth$ and distinct elements $x_\alpha, y_\alpha \in \roth$ for $\alpha<\fc$ such that 
\be 
\itm the sets $F_\alpha, \tilde{F}_\alpha$ are closed under finite intersections;
\itm $\Un_{\beta<\alpha} F_\beta \sub F_\alpha, \Un_{\beta<\alpha} \tilde{F}_\beta \sub \tilde{F}_\alpha$;
\itm $\card{F_\alpha}<\fc, \smallcard{\tilde{F}_\alpha}<\fc$;
\itm $\set{d_\beta, x_\beta}{\beta<\alpha} \le_{F_\alpha} x_\alpha$ and $\set{d_\beta, y_\beta}{\beta<\alpha} {\le_{\tilde{F}_\alpha} }y_\alpha$;
\itm $d_\alpha \le^\ast \max \{x_\alpha\comp, y_\alpha\comp\}$.
\ee 
Let $X_\alpha \coloneqq \set{x_\gamma}{\gamma<\alpha}$; $Y_\alpha \coloneqq \set{y_\gamma}{\gamma<\alpha}$. 
We say that $\alpha$ is $X$-good if $\cU_{\alpha, n}$ is an $\w$-fat cover of $X_\alpha$ for each $n \in \w$. 
We define the notion of a $Y$-good ordinal in the same way. By Lemma~\ref{lem:sbfombfom} for each $X$-good $\alpha$
there are sets $U_{\alpha, n}^X \in \cU_{\alpha,n}, n \in \w$ such that $\set{U_{\alpha,n}^X}{n \in \w}$ is an $\w$-fat cover of {$X_\alpha$.}
The same is true for $Y$-good ordinals.
We want the sequence $\seq{x_\alpha}{\alpha<\mathfrak{c}}$ and the sequence $\seq{y_\alpha}{\alpha<\mathfrak{c}}$ to satisfy the following conditions:
\be 
\itm[(6)] for each $X$-good $\beta < \alpha$ the set $\set{U_{\beta, n}^X}{n \in \w}$ is an $\w$-fat cover of $X_{\alpha}$,
\itm[(7)] for each $Y$-good $\beta < \alpha$ the set $\set{U_{\beta, n}^Y}{n \in \w}$ is an $\w$-fat cover of $Y_{\alpha}$.
\ee 

We begin by setting $F_0=\tilde{F}_0=\{\w\}$ and $X_0=Y_0=\emptyset$.
For $A\in\{X,Y\}$, if $0$ is $A$-good, let
\[
G_{0,m,\emptyset}^A \coloneqq \Un \sset{U_{0,n}^A}{n \in \w \text{ and for each } O \in \cF_m, U_{0, n}^A \cap O \notin \cM} \text{ for } m \in \w
\] and
\[
L_0^A\coloneqq
\bigcap_{m\in\w}G_{0,m,\emptyset}^A;
\]
otherwise, let $L_0^A\coloneqq\roth$.
By Corollary~\ref{cor:omega_fat_corrolary1}(1), the sets $L_0^X$ and $L_0^Y$ are comeager.
Thus
\[
L_0\coloneqq G_0\cap L_0^X\cap L_0^Y
\]
is comeager.
Applying Lemma~\ref{lem:comb_luzin_fin_mod} with $Z=\emptyset$ and with a family of meager sets whose union is $\roth\sm L_0$, we obtain
$x_0,y_0\in\roth$ such that
\[
x_0,y_0\in L_0\cap O_0
\quad\text{and}\quad
d_0\les\max\{x_0\comp,y_0\comp\}.
\]
By Corollary~\ref{cor:omega_fat_corrolary1}(2), conditions~(6) and~(7) are satisfied for $\alpha=1$.
 
Now let $\alpha > 0$ and suppose we have $x_\beta, y_\beta, F_\beta, \tilde{F}_\beta$ such that the conditions
(1)-(7) are satisfied for all $\beta < \alpha$.
Let 
\[
 F=\Un_{\beta<\alpha} F_\beta, \text{ let } \tilde{F}=\Un_{\beta<\alpha} \tilde{F}_\beta, 
\]
let $S \coloneqq \maxfin \set{d_\beta, x_\beta}{\beta<\alpha}$, and let $\tilde{S} \coloneqq \maxfin \set{d_\beta, y_\beta}{\beta < \alpha}$.

Let $Z_{\alpha} \coloneqq \set{s \circ f}{s \in S \cup \tilde{S}, f \in F \cup \tilde{F}}$.
Since $\cov(\cM)$ is regular, it follows that \[\smallcard{S}, \smallcard{\tilde{S}}, \smallcard{F}, \smallcard{\tilde{F}}< \cov(\cM).\]
Thus, $\card{Z_\alpha}< \cov(\cM) = \fc$.

For $A \in \{X,Y\}$, every finite set $F \sub A_{\alpha}$, every $A$-good $\beta \le \alpha$, and every $m \in \w$ define 
\[
 G_{\beta, m, F}^A \coloneqq 
 \Un \set{U_{\beta,n}^A}{
 n\in\w,\ 
 F \sub U_{\beta,n}^A \text{ and for each } O \in \cF_m,\ U_{\beta,n}^A \cap O \notin \cM
 }.
\]

By Corollary~\ref{cor:omega_fat_corrolary1}(1) and the inductive hypothesis, $G_{\beta, m,F}^A$ is comeager for 
each $m \in \w$, each finite set $F \sub A_{\alpha}$, each $A$-good $\beta \le \alpha$, and $A=X,Y$.
Then let \[L' = \bigcap_{\beta<\alpha} G_\beta \cap \bigcap_{\beta <\alpha} \bigl(\roth \sm \{x_\beta, y_\beta\}\bigr)\] and let
\[
	L=L'\cap \bigcap \set{G_{\beta, m, F}^A}{A=X,Y, A\text{-good } \beta \le \alpha, m \in \w, \text{ finite set } F \sub A_{\alpha}}.
\]
Observe that $L=\roth \sm \Un \cM_{\alpha}$ for some family of meager sets $\cM_{\alpha}$ with 
$\card{\cM_{\alpha}}<\fc=\cov(\cM)$. Let $k=\alpha \mod \w$. 
Let $x,y\in\roth$ be functions obtained from Lemma~\ref{lem:comb_luzin_fin_mod} applied to $Z_\alpha$ and $\cM_\alpha$. 
Then there are finite modifications $x_\alpha$ and $y_\alpha$ of $x$ and $y$, respectively such that 
$x_{\alpha}, y_{\alpha} \in L \cap O_k$, 
$d_{\alpha} \le^\ast \max\{x_{\alpha}\comp,y_{\alpha}\comp\}$, and $Z_{\alpha} \le^{\infty} x_{\alpha}, y_{\alpha}$.
By Corollary~\ref{cor:omega_fat_corrolary1}~(2) conditions~(6) and~(7) are preserved. Now we need to show that the remaining conditions are preserved as well. 

The intersection of finitely many elements of the set $F \cup \set{[s \le x_{\alpha}]}{s \in S}$ is infinite:
Since the set $F$ is closed under finite intersections, it is enough to show that for an element $f \in F$ and a finite subset
$S' \sub S$, the intersection $f \cap \bigcap_{s \in S'} [s \le x_{\alpha}]$ is infinite.
Let $f \in F$ and $S' \sub S$ be a finite set. Then
\[
 \bigcap_{s \in S'} [s \le x_{\alpha}]=\bigcap_{s \in S'} \set{n}{s(n) \le x_{\alpha}(n)} \supseteq \smallmedset{n}{\max_{s \in S'}s(n) \le x_\alpha(n)}=[\max[S'] \le x_\alpha].
\]
From the definition of $S$ we have $\max[S'] \in S$. Thus, $\max[S'] \circ f \le^\infty x_\alpha$ and 
hence
\[ 
 \max[S'](f(n)) \le x_{\alpha}(n) \le x_{\alpha}(f(n)),
\]
for infinitely many natural $n$. It follows that the intersection $f \cap [\max[S'] \le x_{\alpha}]$ is infinite and $f \cap \bigcap_{s \in S'} [s \le x_\alpha]$ is infinite as well.
Let $F_\alpha$ be the closure of the set $F \cup \set{[s \le x_\alpha]}{s \in S}$ under finite intersections. By the definition of $S$ we have $d_\beta, x_\beta \in S$ for all ordinals $\beta<\alpha$.
It follows that $\set{d_\beta, x_\beta}{\beta<\alpha} \le_{F_\alpha} x_\alpha$. Analogously, we define the set $\tilde{F}_\alpha$ and show that $\set{d_\beta, y_\beta}{\beta<\alpha} \le_{\tilde{F}_\alpha} y_\alpha$. 
By construction we have $\smallcard{F_{\alpha}}, \smallcard{\tilde{F}_\alpha}<\fc=\cov(\cM)$.
Thus, all conditions are satisfied.

Let $U$ and $\tilde{U}$ be free ultrafilters containing the sets $\Un_{\alpha<\fc}F_\alpha$ and $\Un_{\alpha<\fc}\tilde{F}_\alpha$,
respectively. Define $X \coloneqq \set{x_\alpha}{\alpha < \fc}$ and $Y \coloneqq \set{y_\alpha}{\alpha<\fc}$.
Then the set $X$ is $\sbomfomf$. Take any sequence $\la \cU_{\alpha, n} \colon n \in \w \ra$ of $\w$-fat Borel covers of $X$.
Then $\sset{U_{\alpha, n}^X}{n \in \w}$ is an $\w$-fat cover of $X$. Indeed,
take any finite family of nonempty basic open sets $\cF$ and any finite set $F \sub X$. Then there 
exists $\alpha<\alpha'<\fc$ such that $F \sub X_{\alpha'}$ and 
$\alpha$ is $X$-good. Then the set $\set{U_{\alpha, n}^X}{n \in \w}$ is an $\w$-fat cover of $X_{\alpha'}$.
Thus, there exists $n \in \w$ such that $F \sub U_{\alpha, n}^X$ and $U_{\alpha,n}^X \cap O$ is not meager for 
each $O \in \cF$. Thus, $\set{U_{\alpha, n}^X}{n \in \w}$ is an $\w$-fat cover of $X$. It follows that 
$X$ is $\sbomfomf$.
Similarly, we show that $Y$ is $\sbomfomf$.
Observe that $X \cap O \cap G \neq \emptyset$ for each nonempty basic open set $O$ and 
for each dense $G_\delta$-set in $\roth$. 
Indeed, take any nonempty basic open set $O$ and dense $G_\delta$-set 
$G$. Then $O=O_k$ and $G=G_\alpha$ for some $k \in \w$ and $\alpha<\fc$. Observe that there exists $\alpha' \in (\alpha, \fc)$
such that $k = \alpha' \mod \w$. Then $x_{\alpha'} \in G_\alpha \cap O_k$.
Similarly, one can show that $Y \cap O \cap G \neq \emptyset$ for each nonempty basic open set $O$ and for each dense 
$G_\delta$-set $G$. It follows that both $X \cap O$ and $Y \cap O$ are nonmeager for each nonempty basic open set $O$.
By Lemma~\ref{lem:all_borel_is_fat} $X$ and $Y$ are $\sbomom$ and so are $X \cup \Fin$ and $Y \cup \Fin$.

The set $X$ is a $U$-scale. By construction, $\card{X}=\fc=\fb(U)$.
Fix an element $b \in \roth$. 
We will show that there exist a function $c$ and a set $S \sub X$ of cardinality smaller than $\fb(U)=\fc$ such that 
\[
b \le_{U^+} c \le_U x
\]
for all functions $x \in X \sm S$.
Since $U$ is an ultrafilter, it follows that 
\[U^+ = \sset{a \in \roth}{a\comp \notin U} = \sset{a \in \roth}{a \in U} = U.\]
Hence, it suffices to show that there exists a set $S \sub X$ of cardinality smaller than $\fc$ such that $b \le_U x$ for all $x \in X \sm S$.
 There is an ordinal $\beta<\fc$ such that $b \le^\ast d_\beta$.
For all ordinals $\alpha \in (\beta, \fc)$ we have $d_\beta \le_U x_\alpha$.
Since $U$ contains all cofinite subsets of $\w$, it follows that $b \le_U d_\beta$.
Thus, $b \le_U x_\alpha$ for all $\alpha\in(\beta, \fc)$.
Analogously, $Y$ is a $\tilde{U}$-scale.

It remains to show that $(X \cup \Fin) \x (Y\cup\Fin)$ is not Menger. Observe that the map
$\varphi \colon (X \cup \Fin) \x (Y\cup\Fin) \to \w^\w$ defined by $\varphi(x,y)=\max\{x\comp, y\comp\}$
is continuous. However, by our construction the set
$\set{\max\{x\comp, y\comp\}}{x \in X, y \in Y}$ is dominating. Thus, $(X \cup \Fin) \x (Y\cup\Fin)$ is not Menger.
\epf
Now we present an example showing that even for a
$\fb$-scale we cannot replace the assumption in Theorem~\ref{thm:kappa_fin_unb_x_sbgg_is_sgom}.
The example uses ideas from the paper by Repov\v{s} and Zdomskyy~\cite{RZ}. 

 \bthm\label{thm:counter2}
 It is consistent with CH that there exist a $\fb$-scale set and a set $Y$ satisfying $\sbomom$ whose product 
 space is not Menger.
 \ethm
 To obtain this result, we require the following two straightforward lemmas.

\blem\label{lem:rel_bet_al_incl_and_al_greater}
Let $a,b \in \roth$ be elements such that $\card{a \sm b}<\w$ and $\card{b \sm a} = \w$. Then $b \le^\ast a$.
\elem

Let $a, b \in \roth$.
We say that $a \sub^\ast b$ if $\card{a \sm b}<\w$.
We say that $X = \set{x_\alpha}{\alpha<\fb} \sub \roth$ is a \emph{$\fb$-unbounded tower} if $X$ is unbounded with respect
to $\le^\ast$ and $x_\alpha \sub^\ast x_\beta$ {as well as $x_\beta \not \sub^\ast x_\alpha$} for all $\beta < \alpha < \fb$. 
\blem\label{lem:unbdd_tower_con_a_b_scale}
Every $\fb$-unbounded tower is a $\fb$-scale.
\elem
\bpf[{Proof of Theorem~\ref{thm:counter2}}]
Suppose that $\CH$ holds in the ground model $V$ and fix 
\[
Y = \set{y_\alpha}{\alpha<\w_1} \sub \roth
\]
such that 
\be
\item $y_\beta \sub^\ast y_\alpha$ for all $\beta > \alpha$;
\item $y_{\alpha+1} \sub y_\alpha$ for all $\alpha<\w_1$;
\item $y_\alpha \sm y_{\alpha+1}$ is infinite for all 
$\alpha<\w_1$;
\item For every $y \in \roth$, there exists $\alpha$ with 
$y_\alpha \not \le^\ast y$.
\ee
Then $Y$ is a $\fb$-unbounded tower.
In what follows, we shall work in $V[G]$ where $G$ is 
$\Fn(\w_1, 2)$-generic over $V$. Here $\Fn(\w_1, 2)$ is the standard poset for adding $\w_1$ Cohen reals over $V$.
The set $Y$ is unbounded in $V[G]$ since Cohen reals preserve the unboundedness of ground-model unbounded sets.
Fix an enumeration $\roth=\set{z_\alpha}{\alpha<\w_1}$
and for every $\alpha$, pick $x_\alpha \in [y_\alpha]^\w$
such that $x_\alpha \supset y_{\alpha+1}$, $\card{x_\alpha \sm y_{\alpha+1}}=\card{y_\alpha \sm x_\alpha}=\w$, and 
$z_\alpha \le^\ast (y_\alpha \sm x_\alpha)$. Then $X = \set{x_\alpha}{\alpha<\w_1}$ is a $\fb$-unbounded tower.
By Lemma~\ref{lem:unbdd_tower_con_a_b_scale} $X$ is a $\fb$-scale.
By~\cite[Observation 3.1]{RZ} $(X \cup \Fin) \x Y$ is not Menger.
It therefore suffices to show that $Y$ is $\sbomom$ in $V[G]$. By~\cite[Theorem 18]{ST2002} it suffices to show that $Y^k$ is $\sone(\cB, \cB)$ for all
$k \in \w$.
By~\cite[Theorem 14]{ST2002} it suffices to show that each Borel image of $Y^k$ is Rothberger for all $k \in \w$.

Let $k \in \w$ and let $W$ be an image of $Y^k$ under a Borel map $g$.
Fix a sequence of open covers $\cU_0, \cU_1, \dotsc$ of $W$ 
in $V[G]$, where $\cU_n = \set{U_{n,m}}{m \in \w}$. 
This sequence and a Borel code for $g$ belong to the model $V[G_{\alpha}]$ for some $\alpha < \w_1$ that arises by adding
$\alpha$ Cohen reals to the ground model.
Fix $w \in W$.
Define 
$f_w(n) := \min\set{m \in \w}{w \in U_{n,m}}$ for all $n$.
Since the Borel code for $g$ is in $V[G_\alpha]$, we have $f_w\in V[G_\alpha]$.
If $c \in V[G]$ is any Cohen real over $V[G_\alpha]$, then for
each $w \in W$ the set $\set{n \in \w}{f_w(n)=c(n)}$ is infinite. Thus, the set 
$\set{U_{n, c(n)}}{n \in \w}$ witnesses the Rothberger property of $W$ in $V[G]$.
\epf


\section{Products of $\fd$-concentrated sets in the Miller model}\label{sec:MillerProducts}
As mentioned in the introduction, by  \cite[Theorem~4.7]{SPMProd} in the Miller model,  the product space of finitely many  sets that are $\fd$-concentrated on $\Fin$ is Menger~\cite[Theorem 1.1]{Zdo18}. However, it is consistent with $\ZFC$ that there are sets that are $\fd$-concentrated on $\Fin$ whose product is not Menger~\cite[Theorem 2.7]{ST}.
Let $U$ and $\tilde{U}$ be ultrafilters. We say that 
 $U$ is \emph{near coherent} to $\tilde{U}$ if there is a function $h \in \roth$ such that for each set 
 $u \in U$, we have
 \[
 \bigcup \set{[h(n), h(n+1))}{u \cap [h(n), h(n+1)) \neq \emptyset} \in \tilde{U}.
 \]
The near coherence of filters ($\NCF$), introduced by Blass, is the statement that any two ultrafilters are near coherent~\cite{Blass1986}. This statement is independent of $\ZFC$; for a proof of consistency see~\cite{BS1989}. Consistency of $\neg\NCF$ can be inferred from a result of Blass~\cite[Corollary~15]{Blass1986} and the fact that $\fb=\fd$ if $\CH$ holds.
Let $U \sub \roth$ be an ultrafilter. We say that $B \sub U$
is a base of $U$ if for each 
$u \in U$, there exists $b \in B$
such that $b \sub u$. By 
$\fu$ we denote the minimal
cardinality of an ultrafilter base.
It is a known fact that, if $\NCF$ holds, then $\fu < \fd=\mathfrak b(U)$ for every
ultrafilter $U$; see \cite[Corollary 15, Theorem 16]{Blass1986}.
In particular, $\fd$ is regular under $\NCF$.
In fact, $\NCF$ is equivalent to the statement ``$\fu < \mathfrak b(U)$
for every
ultrafilter $U$'' by \cite[Proposition~2.1(c), Proposition~3.2(e)]{Mil01}.
Note that the set $X$ is $\fd$-unbounded if and only if the set $X \cup \Fin$ is $\fd$-concentrated on $\Fin$~\cite[Lemma 2.3]{ST}.
Recall that a set $X\sub\roth$ is a $U$-scale~\cite[Definition 4.1]{ST} if $\card{X}\geq \bof(U)$ and for any function $d\in\roth$ there is a set $S\sub X$ of cardinality smaller than $\bof(U)$ such that $d \leq_U x$ for all functions $x\in X\sm S$.

\blem \label{lem:dunbdd}
Let $U$ be an ultrafilter such that $\bof(U)=\fd$ and $U$ has a basis of cardinality less than $\fd$.
Each $\fd$-unbounded set in $\roth$ is a $U$-scale.
\elem

\bpf
Suppose that $X\sub\roth$ is a $\fd$-unbounded set that is not a $U$-scale.
Then there is a function $g\in\roth$ such that
\[
\card{\sset{x\in X}{x<_U g}}\geq \fd.
\]
Since every subset of $X$ of cardinality at least $\fd$ is $\fd$-unbounded and $\fd$ is regular, we may assume that $X$ is $\le_U$-bounded by a function $g$.
Let $B\sub\roth$ be a basis of $U$ of size less than $\fd$.
Define $g_b:=\sset{g(n)}{n\in b}$ for all $b\in B$.
Since $X$ is $\fd$-unbounded and $\card{B}<\fd$, there is a function $x\in X$ with $\sset{g_b}{b\in B}<^\infty x$.
Since $x\le_U g$, there is a set $b\in B$ with $b\sub [x\leq g]$. Take a natural number $k$ with $g_b(k)<x(k)$ {and let $n = b(k)$.}
 {Then $n \geq k$ and we have}
\[
g(n)= {g(b(k))} =g_b(k)<x(k)\leq x(n)\leq g(n),
\]
a contradiction.
\epf

\bthm\label{thm:fd_unbdd_in_Miller}
Assume that $\NCF$ holds.
Let $U\sub\roth$ be an ultrafilter.
Each $\fd$-unbounded set in $\roth$ is a $U$-scale.
\ethm

\bpf
By $\NCF$ there exists an ultrafilter $G$ generated by a base $B\sub G$ of cardinality less than $\fd$, and a monotone surjection $\varphi\colon\w\to\w$ such that $F:=\varphi[G]=\varphi[U]$.
Then $\varphi[B]$ is a base for $F$.  Let $X\sub\roth$ be a $\fd$-unbounded set.
By Lemma~\ref{lem:dunbdd}, the set $X$ is an $F$-scale.

Now suppose that, contrary to our claim,  the set $X$ is not a $U$-scale,
i.e., the set $Y=\sset{x\in X}{x\leq_U b}$ has size $\fd$ for some $b\in\roth$.
Set $[k_i,k_{i+1})=\varphi^{-1}(i)$ and
$c(i)=b(k_{i+1})$ for all $i\in\w$. 
Fix $x\in Y$.
Note that if $x(n)\leq b(n)$ for some $n\in [k_i,k_{i+1})$,
then
\[
x(i)\leq x(k_i)\leq x(n)\leq b(n)\leq b(k_{i+1})=c(i)\quad\text{and}\quad \varphi(n)=i.
\]
It follows that $\varphi\bigl[[x\leq b]\bigl]\sub [x\leq c]$.
Since $[x\leq c]\in U$, we have $\varphi\bigl[[x\leq b]\bigl]\in F$, and thus $x\leq_F c$.
We conclude that $Y$ is bounded with respect $\leq_F$, a contradiction.
\epf

In the Miller model, by~\cite[Theorem~2]{BL1989} and results of~\cite{BS1989}, the inequality $\fu < \fg$ holds,  see~\cite{BL1989} for the definition of $\fg$.
As a result,  the NCF holds in the Miller model,
as it was established in~\cite{L1992}.

\bcor \label{cor:6.3}
In the Miller model, for every positive natural number $k$ and $\fd$-concentrated sets $Y_i\sub\PN$ on  
$\fin$ for $i<k$, the product $\prod_{i
< k} Y_i$ satisfies $\sgw$.
\ecor

\bpf
Let $X_i=Y_i\sm\fin$ and note that each set $X_i$ is $\fd$-unbounded, where $i<k$.
Thus, the set $X=\bigcup_{i<k} X_i$ is also $\fd$-unbounded, and hence by Theorem~\ref{thm:fd_unbdd_in_Miller} it is a $U$-scale for an ultrafilter $U$. 
Now, Proposition~\ref{lem:F-scale_is_b(f)_unbdd}
together with $\bof(U)=\fd$ implies that $X$ is  a
$\fd$-fin-unbounded set.
Finally, 
Proposition~\ref{prp:mainmich} ensures that 
$\prod_{i
< k} Y_i$ satisfies $\sgw$.
\epf

Let us note that it is not essential in Corollary~\ref{cor:6.3} on which 
countable subsets of the Cantor set the sets $Y_i$
are concentrated since the Cantor set is countable dense homogeneous, i.e., 
for any two countable dense subsets of $\PN$ there exists an autohomeomorphism of $\PN$ mapping one of them onto the other one (see~\cite{ArhMil14} and references therein).

\section{Comments and open problems}

\subsection{Characterization of $\sone(\Gamma, \Lambda), \sone(\Gamma, \Omega)$, and $\sone(\Gamma, \Gamma)$}
In this section, we present the characterizations of $\sone(\Gamma, \Lambda), \sone(\Gamma, \Omega)$, and 
$\sone(\Gamma, \Gamma)$. The characterization of $\sgg$ appears implicitly in the paper by Peng~\cite{P2023}.
Here we state his result more explicitly and devise characterizations of $\sone(\Gamma, \Lambda)$ and $\sone(\Gamma, \Omega)$
using his ideas. We denote by $\w+1$ the space $\w \cup \{\w\}$ endowed with the order topology.
Consequently, $\w+1$ is a one point compactification of $\w$.
\blem\cite[Proof of Proposition~6.]{P2023}\label{lem:Continuous_function_from_X_to_omega+1}
Let $X$ be a zero-dimensional, separable metric space and $\eseq{\cU}\in\Gamma(X)$, where $\cU_n = \sset{U_{n,m}}{m \in \omega}$ for all $n$.
Then for all natural numbers $n$ and $m$ there is a continuous function $f_{n,m} \colon X \to \w+1$ such that 
$X \sm U_{n,m}=f_{n,m}^{-1}[\{\w\}]$.
\elem

\bpf
Fix natural numbers $n$ and $m$.
Then $U_{n,m} = \Un_{i \in \omega} A_{n,m,i}$, where each $A_{n,m, i}$ is 
clopen and 
$A_{n,m,i} \sub A_{n,m,i+1}$ for all
$i \in \w$.
Define $f_{n,m} \colon X \to \omega+1$ by 
\beqs
f_{n,m}(x) \coloneqq \begin{cases}
 i,&\text{ if } x \in A_{n,m,i} \sm \Un_{j<i} A_{n,m,j},\\ 
 \w,& \text{ if } x \not\in U_{n,m}\\
\end{cases}.
\eeqs

Fix a nonempty basic open set $U \sub \w+1$. If $\w \notin U$, then 
\[f_{n,m}^{-1}[U] = \Un_{i \in U} f_{n,m}^{-1}[\{i\} ]= \Un_{i \in U} \biggl(A_{n,m,i} \sm \Un_{j<i} A_{n,m,j}\biggr)\]
is an open subset of $X$.
Assume that $\w \in U$.
Then $U = (\w+1) \sm F$ for some finite subset $F \sub \w$ and  
\[
f_{n,m}^{-1}[U] = X \sm \Un_{i \in F} f_{n,m}^{-1}[\{i\}] = X \sm \Un_{i \in F}\biggl(A_{n,m,i} \sm \Un_{j<i} A_{n,m,j}\biggr).
\] 
Since $\Un_{i \in F} \bigl(A_{n,m,i} \sm \Un_{j<i} A_{n,m,j}\bigr)$ is clopen, it follows
that $f_{n,m}^{-1}[U]$ is open.

\epf
Fix $a\in\NN$.
We identify $a$ with its graph, regarded as a subset of $\w\x\w$.
We denote by $\mathrm{cF}(a)$ the filter consisting of all 
cofinite subsets of $a$.
We denote by $\cX$ the set of all functions  $f\in (\w+1)^{\w\x\w}$ such that the sets $\sset{m}{f(n,m)=\w}$ are finite for all natural numbers $n$.
Consider $\cX$ with a subspace topology of the space $(\w+1)^{\w\x\w}$ with the product topology.
\bthm\label{thm:Characterization}
Let $X$ be a zero-dimensional, separable metric space.
\be
\item The space $X$ is $\sgg$ if and only if for each continuous 
$\varphi \colon X \to \cX$ there is $a \in {\w^\w}$ such that $\sset{a \sm \varphi(x)^{-1}\{\omega\}}{x \in X} \sub \mathrm{cF}(a)$;
\item The space $X$ is $\sgom$ if and only if for each continuous 
$\varphi \colon X \to \cX$ there is $a \in {\w^\w}$ such that $\sset{a \sm \Un_{x \in F}\varphi(x)^{-1}\{\omega\}}{F \in [X]^{<\omega}} \sub [a]^{\omega}$;
\item The space $X$ is $\sone(\Gamma,\Lambda)$ if and only if for each continuous 
$\varphi \colon X \to \cX$ there is $a \in {\w^\w}$ such that 
 $\sset{a \sm \varphi(x)^{-1}\{\omega\}}{x \in X}\sub [a]^{\omega}$.
\ee
\ethm

\bpf
($\Rightarrow$)
Let $\cA\in\{\Gamma,\Omega,\Lambda\}$ and assume that $X$ is $\sone(\Gamma, \cA)$.
Fix a continuous function $\varphi \colon X \to \cX$. 
Let 
\[
U_{n,m}\coloneqq \sset{g \in \cX}{g(n,m)<\omega}
\]
for all natural numbers $n, m$.
Let
\[
T\coloneqq \sset{n\in\w}{\varphi[X]\sub U_{n,m}\text{ for some }m\in\w}.
\]
For every $n\in T$, choose $a(n)\in\w$ such that $\varphi[X]\sub U_{n,a(n)}$.
For every $n\in\w\sm T$, the family
\[
\sset{U_{n,m}\cap\varphi[X]}{m\in\w}
\]
is a $\gamma$-cover of $\varphi[X]$. Indeed, since $\varphi[X]\sub\cX$, for every $g\in\varphi[X]$ we have
$g\in U_{n,m}$ for all but finitely many $m$, and, by the definition of $T$,
$U_{n,m}\cap\varphi[X]\neq\varphi[X]$ for every $m\in\w$.

If $\w\sm T$ is infinite, enumerate it as $\w\sm T=\sset{n_k}{k\in\w}$.
Since the set $\varphi[X]$ is $\sone(\Gamma,\cA)$, there are $a(n_k)\in\w$ for all $k\in\w$ such that
\[
\sset{U_{n_k,a(n_k)}\cap\varphi[X]}{k\in\w}\in\cA(\varphi[X]).
\]
If $\w\sm T$ is finite, define $a(n)$ arbitrarily for all $n\in\w\sm T$.
Thus, $a\in\w^\w$ is defined.

(1) Assume that $\cA=\Gamma$.
For each $x\in X$, we have $\varphi(x)\in U_{n,a(n)}$ for every $n\in T$ and for all but finitely many $n\in\w\sm T$.
Consequently,
\[
\sset{a \sm \varphi(x)^{-1}\{\w\}}{x \in X} \sub \mathrm{cF}(a).
\]

(2) Assume that $\cA=\Omega$.
Let $F\in[X]^{<\omega}$.
If $T$ is infinite, then
\[
\sset{(n,a(n))}{n\in T}
\sub
a \sm \Un_{x \in F}\varphi(x)^{-1}\{\w\}.
\]
If $T$ is finite, then $\w\sm T$ is infinite and
\[
\sset{U_{n,a(n)}\cap\varphi[X]}{n\in\w\sm T}\in\Omega(\varphi[X]).
\]
Since every finite subset of a space is contained in infinitely many members of an $\omega$-cover, there are infinitely many $n\in\w\sm T$ such that $\varphi(x)\in U_{n,a(n)}$ for all $x\in F$.
In either case,
\[
\sset{a \sm \Un_{x \in F}\varphi(x)^{-1}\{\w\}}{F \in [X]^{<\omega}} \sub [a]^{\omega}.
\]

(3) Assume that $\cA=\Lambda$.
Fix $x\in X$.
If $T$ is infinite, then
\[
\sset{(n,a(n))}{n\in T}
\sub
a\sm\varphi(x)^{-1}\{\w\}.
\]
If $T$ is finite, then $\w\sm T$ is infinite and
\[
\sset{U_{n,a(n)}\cap\varphi[X]}{n\in\w\sm T}\in\Lambda(\varphi[X]),
\]
and therefore $\varphi(x)\in U_{n,a(n)}$ for infinitely many $n\in\w\sm T$.
In either case,
\[
\sset{a \sm \varphi(x)^{-1}\{\w\}}{x \in X}\sub [a]^{\omega}.
\]

($\Leftarrow$)
Let $\eseq{\cU}\in\Gamma(X)$, where $\cU_n=\sset{U_{n,m}}{m \in \omega}$ for all $n$.
By Lemma~\ref{lem:Continuous_function_from_X_to_omega+1} for all natural numbers $n$ and $m$ there is a continuous function
$f_{n,m} \colon X \to \w+1$ such that $X \sm U_{n,m}=f_{n,m}^{-1}[\{\w\}]$.
Then the map $\varphi \colon X \to \cX$, defined by $\varphi(x)(n,m)\coloneqq f_{n,m}(x)$, is continuous.

(1) There is $a \in {\w^\w}$ such that
\[
\sset{a \sm \varphi(x)^{-1}\{\omega\}}{x \in X} \sub \mathrm{cF}(a).
\]
Fix $x\in X$.
We have $x \in U_{n, a(n)}$ for all but finitely many $n$. Thus,
\[
\sset{U_{n, a(n)}}{n \in \omega}\in\Gamma(X).
\]

(2) There is $a \in {\w^\w}$ such that
\[
\sset{a \sm \Un_{x \in F}\varphi(x)^{-1}\{\w\}}{F \in [X]^{<\omega}} \sub [a]^{\omega}.
\]
Let $F$ be a finite subset of $X$.
Then $a \sm \Un_{x \in F} \varphi(x)^{-1}\{\omega\}$ is an infinite subset of $a$.
Thus, there is $n$ such that $\varphi(x)(n, a(n))<\omega$ for all $x \in F$, and hence $F \sub U_{n,a(n)}$.
It follows that 
\[
\sset{U_{n, a(n)}}{n \in \omega}\in\Om(X).
\]

(3) There is $a \in {\w^\w}$ such that
\[
\sset{a \sm \varphi(x)^{-1}\{\omega\}}{x \in X} \sub [a]^{\omega}.
\]
Fix $x \in X$.
Then $a \sm \varphi(x)^{-1}\{\omega\}$ is an infinite subset of $a$, and thus there are infinitely many $n$ such that $\varphi(x)(n, a(n))<\omega$.
We have $x \in U_{n,a(n)}$ for infinitely many $n$.
It follows that 
\[
\sset{U_{n, a(n)}}{n \in \omega}\in\Lambda(X).\qedhere
\]
\epf
\subsection{Properties $\sgom$ and $\sgo$}

The following fact is an easy and useful observation about spaces all of whose finite powers satisfy $\sgw$.
Surprisingly, it is not mentioned in the literature on combinatorial covering properties.

\bprp\label{prp:sgosgw}
Let $X$ be a space. The space $X^k$ satisfies $\sgo$ for each natural number $k \ge 1$ if and only if $X^k$ satisfies $\sgw$ for each natural number $k \ge 1$.
\eprp
\bpf
It is enough to prove that if a space $X^k$ is $\sgo$ for each natural number $k \ge 1$, then it is $\sgom$.
Let $\eseq{\cU}\in \Gamma(X)$. Let us reindex this sequence as $\cU_{k,m}$ where $k \ge 1$ and $m$ are natural numbers.
Let $\cU_{k,m}':=\sset{U^k}{U \in \cU_{k,m}}$. Then for all natural numbers $k \ge 1$ and $m$, $\cU_{k,m}'\in \Ga(X^k)$.
Since $X^k$ is $\sgo$ for each natural $k \ge 1$, it follows that for each natural number $k \ge 1$, there exist sets $U_{k,m} \in \cU_{k,m}, m \in \w$ such that 
$\set{U_{k,m}^k}{m \in \w}$ is a cover of $X^k$. 

Let $\cU':=\set{U_{k,m}}{k,m \in \w}$. Then $\cU'$ is an $\omega$-cover of $X$. 
Indeed, take any finite $F=\{x_0, x_1, \dotsc, x_{k-1}\} \sub X$. Then there exists $m \in \w$ such that $(x_0, x_1, \dotsc, x_{k-1}) \in U_{k,m}^{k}$. 
Thus, $F \sub U_{k,m} \in \cU'$.
\epf
This fact allows us to easily prove a more basic version of Theorem~\ref{thm:kappa_fin_unb_x_sbgg_is_sgom}~(2) where we take $Y$ to be a one-point space.

\subsection{Michael topology}

Note that all results in which the space is endowed with the Michael topology hold when we consider any topology in which the neighborhoods of $\Fin$ are the same as in the standard topology on $\PN$. 

Let $C$ be the standard middle-third Cantor set in the real line, and let $C_{\mathrm S}$ denote the set $C$ endowed with the topology induced by the Sorgenfrey line~\cite{Sorgenfrey}.
Using the standard homeomorphism between $\PN$ and $C$, we introduce the Sorgenfrey topology on $\PN$ and denote the resulting space by $\PN_{\mathrm S}$.
Since the points of $\Fin$ correspond to the left endpoints of the intervals appearing in the construction of $C$, their open neighborhoods in $\PN$, $\PNM$, and $\PN_{\mathrm S}$ are the same as in the usual real-line topology.
Moreover, for every $X\sub\roth$, the Michael topology on $X\cup\Fin$ is finer than the Sorgenfrey topology.
Thus, the identity map from $\mich{X\cup\Fin}$ to $(X\cup\Fin)_{\mathrm S}$ is continuous.
Consequently, all results above for the Michael topology transfer to the Sorgenfrey topology, since the covering properties under consideration are preserved by continuous images.
For example, one can prove the following result.
\bthm 
Let $\kappa$ be a regular cardinal number with $\kappa \leq \fd$, let $X$ be a $\kappa$-fin-unbounded set, let $k$ be a natural number, and let $Y$ be a set satisfying $\sbgg$. Then the space $(X \cup \Fin)_{\mathrm{S}}^k \x Y$ satisfies $\sgom$.
\ethm 
\subsection{Open problems}

\bprb \label{7.5prb}
Is it consistent with $\ZFC$ that there is a filter $F$ and $F$-scales $X$ and $Y$ such that the product space $(X\cup\Fin)\x (Y\cup\Fin)$ is not $\sgw$? 
\eprb
By Proposition~\ref{prp:mainmich}, a pair of $F$-scales $X$ and $Y$ such that $X\cup Y$ is $\mathfrak b(F)$-fin-unbounded cannot witness an affirmative answer to Problem~\ref{7.5prb}.

\bprb 
Let $A \sub \roth$ satisfy $\card{A} < \fb$ (or $\card{A}<\fd)$, and let $Y$ have the $\gamma$-property.
Does the product space $A \x Y$ satisfy $\sgo$?
\eprb
\bprb 
Let $X$ be a $\fb$-unbounded tower and $Y$ be a set with the $\gamma$-property.
Does the product space $(X \cup \Fin)\x Y$ satisfy $\sgo$?
\eprb
\bprb 
Let $X$ be a $\fb$-unbounded tower and $Y$ be a set with the $\gamma$-property.
Does the product space $(X \cup \Fin) \x Y$ satisfy $\sgom$?
\eprb
\bprb 
Let $\kappa$ be a regular cardinal number and $X$ be a $\kappa$-fin-unbounded scale set. Let $Y$ be a set with the $\gamma$-property.
Does the product space $(X \cup \Fin)\x Y$ satisfy $\sgom$?
\eprb
Let $\pr(\sgom)$ be the minimal cardinality of a space that is not productively 
$\sgom$ in the class of metric spaces.
\bprb
Let $\kappa$ be a regular cardinal number with $\kappa \le \pr(\sgom)$ and let $X$ be a  $\kappa$-fin-unbounded set. Let $Y$ be a space satisfying the $\gamma$-property.
Does the product space $(X \cup \Fin) \x Y$ satisfy $\sgom$?
\eprb

\section{Acknowledgments}

We are very grateful to the referees for their valuable comments and suggestions, which substantially improved the quality of the paper and the presentation of the material.
\mbox{}\\

The research of the first two authors was funded by the National Science Center, Poland, under the Weave-UNISONO call of the Weave programme, project \enquote{Set-theoretic aspects of topological selections} (2021/03/Y/ST1/00122).
The research of the third author was funded in full by the Austrian Science Fund (FWF) [10.55776/I5930 and 10.55776/PAT5730424]. 
This research was completed while the first author was a doctoral
candidate in the Interdisciplinary Doctoral School at the 
{\L}\'{o}d\'{z} University of Technology, Poland.

\end{document}